%% file: lmsf.tex
\documentclass[a4paper]{article}
\input roby3.tex

\usepackage{makeidx}
\usepackage{amsfonts}
\usepackage{amsmath}
\usepackage{latexsym}
\usepackage{amssymb}
\usepackage[dvips]{graphicx}

\title{Applications of Minor Summation Formula III,
Pl\"ucker Relations, Lattice Paths and Pfaffian Identities}

\author{
Masao ISHIKAWA
\ and\ Masato WAKAYAMA
}

\date{\empty}

\newtheorem{thm}{Theorem}[section]
\newtheorem{defi}[thm]{Definition}

\newtheorem{rem}[thm]{Remark}
\newtheorem{ex}[thm]{Example}
\newtheorem{prop}[thm]{Proposition}
\newtheorem{lem}[thm]{Lemma}
\newtheorem{cor}[thm]{Corollary}

\begin{document}
\maketitle

\begin{center}
${}^*$Department of Mathematics, 
Faculty of Education, 
Tottori University
\medskip

${}^\dagger$Faculty of Mathematics,
Kyushu Unversity
\end{center}
\medskip

\abstract{
The initial purpose of the present paper is to provide a combinatorial
proof of the minor summation formula of Pfaffians in \cite{IW1}
based on the lattice path method. The second aim is to study 
applications of the minor summation formula for 
obtaining several identities. Especially, 
a simple proof of Kawanaka's formula concerning a $q$-series identity
involving Schur functions
\cite{Ka1} and of the identity in \cite{Ka2} which is regarded as a
determinant version of the previous one are given.

}

{
\small
\tableofcontents
}

\input lmsfd.tex

\input lmsf1.tex

\input lmsf2.tex

\input lmsf3.tex

\input lmsf4.tex

\input lmsf6.tex

\input lmsf7.tex

\input appendix.tex

\bigbreak
\noindent
{\bf Note Added in Proof }
In a private communication (``Kawanaka's $q$-Cauchy identity and rational universal character''),
S.~Okada has suggested that \thetag{\ref{eqn:kawanaka}} can be proved
by the rational universal characters of the classical groups
(see K.~Koike's paper:
``On the decomposition of tensor products of the representations of the classical groups: 
By means of the universal characters'',
Adv. Math.
{\bf  74} (1989), 57--86.
).

\input lmsfr.tex
\medskip
\parindent=0mm

Masao ISHIKAWA, Department of Mathematics, Faculty of Education, \\
Tottori University, Tottori 680-8551, Japan

\smallskip

E-mail address: ishikawa@fed.tottori-u.ac.jp
\medskip

Masato Wakayama, Faculty of Mathematics,
Kyushu Unversity, Hakozaki, Fukuoka 812-8581, Japan

\smallskip

E-mail address: wakayama@math.kyushu-u.ac.jp

\end{document}

%% file: roby3.tex
\newdimen\Squaresize \Squaresize=20pt
\newdimen\thickness \thickness=1pt         
                                                    
\def\Square#1{\hbox{\vrule width \thickness
   \vbox to \Squaresize{\hrule height \thickness\vss                            
      \hbox to \Squaresize{\hss#1\hss}
   \vss\hrule height\thickness} 
\unskip\vrule width \thickness} 
\kern-\thickness}                                                            
                               
\def\vsquare#1{\vbox{\Square{$#1$}}\kern-\thickness}
\def\blank{\omit\hskip\Squaresize}

\def\fibyoung#1{\let\\=\cr              
\vbox{\smallskip\offinterlineskip
\halign{&\vsquare{##}\cr #1}}\,}

\def\borderlessrect#1#2{\hbox{\hskip \thickness
   \vbox to \Squaresize{\vskip \thickness \vss
      \hbox to #2 {\hss #1\hss}
   \vss\vskip\thickness} 
\unskip\hskip \thickness} 
\kern-\thickness}                                                            
                               
\def\vborderlessrect#1#2{\vbox{\borderlessrect{$#1$}{#2}}\kern-\thickness}

\def\borderless#1{\omit\vborderlessrect{#1}{\Squaresize}}
\def\borderlessrc#1#2{\omit\vborderlessrect{#1}{#2}}

\def\msquare#1{\vbox{\hbox{\vrule width \thickness
   \vbox to \Squaresize{\hrule height \thickness
      \hbox to \Squaresize{\hfil{\sevenrm #1}}
   \vfil\hrule height\thickness}
\unskip\vrule width \thickness}
\kern-\thickness}\kern-\thickness}

\def\twosquare#1#2{\vbox{\hbox{\vrule width \thickness 
   \vbox to \Squaresize{\hrule height \thickness
      \hbox to \Squaresize{\hfil{\sevenrm #1}}\vss
      \hbox to \Squaresize{\hss{#2}\hss}
   \vfil\hrule height\thickness}
\unskip\vrule width \thickness}
\kern-\thickness}\kern-\thickness}

\def\twoblank#1#2{\vbox{\hbox{
   \vbox to \Squaresize{\vskip 2pt
      \hbox to \Squaresize{\hfil{\sevenrm #1}\ }\vss
      \hbox to \Squaresize{\hss{#2}\hss}
   \vfil}\unskip\kern-\thickness}
}\unskip\kern-\thickness}

\def\young#1{
\def\>{\blank}
\def\<{\borderless}
\def\*{\borderlessrc}
\def\p{\omit\msquare}
\def\t{\omit\twosquare}
\def\b{\omit\twoblank}
\let\\=\cr 
\vbox{\smallskip\offinterlineskip
\halign{&\vsquare{##}\cr #1}}}

\newdimen\smsquaresize \smsquaresize=12pt
\newdimen\smthickness \smthickness=.5pt
\font\smcellfont=cmss8 scaled \magstep0

\def\smsquare#1{\hbox{\vrule width \smthickness
   \unskip\vbox to \smsquaresize{\hrule height \smthickness\vss
      \hbox to \smsquaresize{\hss{\smcellfont #1}\hss}
   \vss\hrule height\smthickness} 
\unskip\vrule width \smthickness} 
\kern-\smthickness}

\def\smvsquare#1{\vbox{\smsquare{$#1$}}\kern-\smthickness}
\def\blank{\omit\hskip\smsquaresize}

\def\smyoung#1{\let\\=\cr 
\vbox{\smallskip\offinterlineskip
\halign{&\smvsquare{##}\cr #1}}}
\newdimen\vsmsquaresize \vsmsquaresize=10pt
\newdimen\vsmthickness \vsmthickness=.5pt
\font\vsmcellfont=cmsl8 scaled \magstep0
\font\vsmletterfont=cmr6 scaled \magstep0

\def\vsmsquare#1{\hbox{\vrule width \vsmthickness
   \unskip\vbox to \vsmsquaresize{\hrule height \vsmthickness\vss
      \hbox to \vsmsquaresize{\hss{\vsmcellfont #1}\hss}
   \vss\hrule height\vsmthickness} 
\unskip\vrule width \vsmthickness} 
\kern-\vsmthickness}
\def\vsmvsquare#1{\vbox{\vsmsquare{#1}}\kern-\vsmthickness}
\def\vsmblank{\omit\hskip\vsmsquaresize}
\def\vsmborderless#1{\hbox{\hskip \vsmthickness\unskip
   \vbox to \vsmsquaresize{\vss
      \hbox to \vsmsquaresize{\hss{\vsmletterfont #1}\hss}
   \vss} 
\unskip\hskip \vsmthickness} 
\kern-\vsmthickness}                                                            \def\vsmvborderless#1{\vbox{\vsmborderless{#1}}\kern-\vsmthickness}

\def\vsmyoung#1{
\def\>{\vsmblank}
\def\<{\omit\vsmvborderless}
\let\\=\cr 
\vbox{\smallskip\offinterlineskip
\halign{&\vsmvsquare{##}\cr #1}}}

%% file: lmsfd.tex



\def\rdots{\mathinner{\mkern1mu\raise1pt\hbox{.}\mkern2mu\raise4pt\hbox{.}\mkern2mu\raise7pt\hbox{.}\mkern1mu}}

\def\defterm#1{{\sl #1}\/}

\def\module{\operatorname{mod}}

\def\ep{\varepsilon}
\def\lam{\lambda}

\def\Comp{\Bbb{C}}
\def\Nat{\Bbb{N}}
\def\Int{\Bbb{Z}}
\def\Pos{\Bbb{P}}

\def\Pf{\operatorname{Pf}}
\def\pf{\operatorname{pf}}
\def\odd{\operatorname{odd}}
\def\even{\operatorname{even}}
\def\sgn{\operatorname{sgn}}
\def\MOD{{\,\operatorname{mod}\,}}

\def\half{\frac12}
\def\trans{{}^t\!}
\def\cc#1{{\ooalign{\hfil\raise-.02ex\hbox{#1}\hfil\crcr\mathhexbox20D}}}
\def\cl{\mathhexbox20D}
\def\cb{\bullet}

\def\diag{\operatorname{diag}}

\def\Ind#1{\Cal{I}_{#1}}





\def\id{\operatorname{id}}
\def\wt{\operatorname{wt}}

          




%% file: lmsf1.tex
%
%
%
%
%
%

\newcount\chapterno \chapterno=1
\newcount\definitionno \definitionno=0
\newcount\propositionno \propositionno=0
\newcount\lemmano \lemmano=0
\newcount\theoremno \theoremno=0
\newcount\corollaryno \corollaryno=0
\newcount\exampleno \exampleno=0
\newcount\remarkno \remarkno=0
\newcount\assumptionno \assumptionno=0
\newcount\claimno \claimno=0
\newcount\conjectureno \conjectureno=0
\newcount\problemno \problemno=0
\newcount\equationno \equationno=0


\section{Introduction}

Recently, applications of the minor summation formula 
presented in \cite{IW1} have been made in several directions, e.g., to study 
a certain 
limit law for shifted Schur measures in \cite{TW1}, 
to find an explicit description of the skew-Capelli 
identity in \cite{KW}, and to generalize further 
the so-called Littlewood formulas, for instance in \cite{IW4}
(see also \cite{La}, \cite{JZ1}, \cite{JZ2}), etc. 
Moreover, the formula has been generalized to 
the case of hyperpfaffians in \cite{LT2} (see also \cite{LT1}).

In this paper, therefore, we treat again the minor summation formulas of Pfaffians 
and derive several basic formulas concerning Pfaffians from certain 
combinatorial theoretical points of view.
In order to develop the study of these formulas nicely, 
we present also a Pfaffian version of 
the Lewis Carroll formula  (Dodgson's identity) and of the Pl\"ucker relations.
The proof and the discussion concerning these formulas have not been 
developed sufficiently in our previous papers because they
are not directly related to the actual proofs in our applications
of the minor summation formulas to obtaining various generating functions of
the Schur functions, etc., whereas they have been studied from the early
stage of the research.

One of the main purpose of the present paper is to prove the minor summation formulas using the lattice path method (see \cite{Ste}) combined with 
the Lewis Carroll formula for Pfaffians.
The proof thus obtained enables us to provide a combinatorial interpretation of the minor summation formula through lattice paths. 
The point is that the Lewis Carroll formula for Pfaffians plays efficiently  
to reduce the proof of the minor summation formulas 
comparing with a similar combinatorial discussion given in \cite{Ste}, 
and the present proof, consequently, may explain 
the meaning of the formulas more clearly. 

The paper is organized as follows.
In Section~2 we present two Pfaffian identities which may be 
called a Pfaffian version of the Lewis Carroll formula 
and the Pl\"ucker relations respectively,
and in Section~3 we formulate the various type of minor summation formulas,
where the Pfaffian 
analogue of the Lewis Carroll formula plays a key role at the derivation.
In fact, we define the notion of the matrix formed by copfaffians 
which may sound abuse of languages,
but we need to define a Pfaffian counterpart of the matrix of cofactors 
in the determinant theory. Then we can 
get a new expression of the minor summation formulas by 
employing these matrices of copfaffians.
Without the notion of the matrices of copfaffians, 
it seems very hard to 
discover Gessel-Viennot type combinatorial proofs developed in Section~4,
which simplify the lattice path proof of the minor summation formulas.

In Section~5 and Section~6 we give certain applications of the minor 
summation formulas to $q$-series. Actually, 
in Section~5 we show that Kawanaka's $q$-Littlewood identity is 
easily derived from the minor summation formulas,
and in Section~6 we show that Kawanaka's $q$-Cauchy identity is proved 
by the Binet-Cauchy formula with some combinatorics.
Furthermore, we shall give  
some variant of the Sundquist formula obtained in \cite{Su2}  
which is considered as a two variable Pfafffian identity.
We put this in the Appendix because the formula is not 
a direct application of the minor summation formula. 

%% file: lmsf2.tex
%
%
%
%
%
%


\section{The Lewis Carroll formula and the Pl\"ucker relations}

We provide a Pfaffian version of Lewis Carroll's formula (Dodgson's identity).
We first recall the so-called Lewis Carroll formula,
or known as Jacobi's formula which is an identity among the minor determinants.
The reader can find a restricted version of this identity and related topics in \cite{B1} and \cite{RR}.
Furthermore, we present the Pl\"ucker's relations.
The latter relations are also treated in \cite{DW},
and in \cite{Kn} they are called the (generalized) basic identity.
We give a brief proof of ordinary Lewis Carroll's formula,
which will be needed to establish the Pfaffian version,
 to make also this paper self-contained.
We only use Cramer's formula to prove it.

Let us denote by $\Nat$ the set of non-negative integers,
and by $\Int$ the set of integers. 
Let $[n]$ denote the set $\{1,2,\dots,n\}$ for a positive integer $n$.
For any finite set $S$ and any nonnegative integer $r$,
let $\binom{S}{r}$ denote the set of all $r$-element subsets of $S$.
For example, $\binom{[n]}{r}$ stands for the set of all multi-indices $\{i_1,\dots,i_r\}$
such that $1\leq i_1<\dots<i_r\leq n$.
Let $n$, $M$ and $N$ be positive integers such that $n\le M,N$
and let $T$ be any $M$ by $N$ matrix.
For any multi-indices
$I=\{i_1,\ldots,,i_n\}\in\binom{[M]}{n}$ and
$J=\{j_1,\dots,j_n\}\in\binom{[N]}{n}$,
let $T^{I}_{J}=T^{i_1\dots i_n}_{j_1\dots j_n}$ be the sub-matrix of $T$
obtained by picking up the rows indexed by $I$ and the columns indexed by $J$, i.e.,
\begin{equation*}
T^{I}_{J}
=\begin{pmatrix}
t_{i_1j_1}&\hdots&t_{i_1j_n}\\
\vdots&\ddots&\vdots\\
t_{i_nj_1}&\hdots&t_{i_nj_n}\\
\end{pmatrix}.
\end{equation*}
In the case of $n=M$ and $I=[M]$,
we omit $I$ from the above expression
and write $T_{J}$ for $T^{I}_{J}$,
when there is no possibility of confusion.
Similarly we may write $T^{I}$ for $T^{I}_{J}$ if n=N and $J=[N]$.

Though the notion of Pfaffians is less familiar than that of determinants,
it is also well-known that the Pfaffian (of a skew-symmetric matrix) is 
expressed  as a square root of the determinant of the corresponding matrix. 
We recall then first a more combinatorial definition of
Pfaffians presented in \cite{Ste}.
Let $\frak S_n$ be the symmetric group 
on the set of the letters $1,2,\ldots,n$, and for each permutation $\sigma\in
\frak S_n$ 
let $\sgn\sigma$ stand for $(-1)^{\ell(\sigma)}$, where
$\ell(\sigma)$ denotes the number of inversions in $\sigma$.

In this paper we use the symbol $\{ i_1, i_2, \dots i_r\}_{<}$ for
the set $\{ i_1, i_2, \dots i_r\}$ with the relation $i_1<i_2<\dots<i_r$.
Let $n=2r$ be an even integer and let 
 $A=(a_{ij})_{1\le i,j\le n}$ be 
an $n$ by $n$ skew symmetric matrix (i.e. $a_{ji}=-a_{ij}$),
whose entries $a_{ij}$ are in a commutative ring.
The \defterm{Pfaffian} $\Pf(A)$ of $A$ is defined by
\begin{equation}
\label{def_pfaffian}
\Pf(A)=\sum \epsilon(\sigma_{1},\sigma_{2},\hdots,\sigma_{n-1},\sigma_{n})\,
a_{\sigma_{1}\sigma_{2}} \dots a_{\sigma_{n-1}\sigma_{n}},
\end{equation}
where the summation is over all partitions $\{\{\sigma_{1},\sigma_{2}\}_{<},\hdots,\{\sigma_{n-1},\sigma_{n}\}_{<}\}$
of $[n]$ into $2$-element blocks,
and $\epsilon(\sigma_{1},\sigma_{2},\hdots,\sigma_{n-1},\sigma_{n})$ denotes the sign of the permutation
\begin{equation*}
\begin{pmatrix}
1&2&\cdots&n-1&n\\
\sigma_{1}&\sigma_{2}&\cdots&\sigma_{n-1}&\sigma_{n}
\end{pmatrix}.
\end{equation*}
For instance, when $n=4$, the equation above reads:
\begin{equation*}
\Pf\begin{pmatrix}
      0& a_{12}& a_{13}& a_{14}\\
-a_{12}&      0& a_{23}& a_{24}\\
-a_{13}&-a_{23}&      0& a_{14}\\
-a_{14}&-a_{24}&-a_{34}&      0
\end{pmatrix}
=a_{12}a_{34}-a_{13}a_{24}+a_{14}a_{23}.
\end{equation*}
Note that a skew symmetric matrix $A=(a_{ij})_{1\leq i,j\leq n}$ is determined by its upper triangular entries $a_{ij}$ for $1\leq i<j\leq n$.

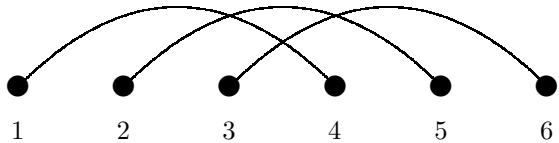
\begin{figure}[b]
\begin{center}
\begin{picture}(200,60)
\put(  0,20){\circle*{8}}
\put( 40,20){\circle*{8}}
\put( 80,20){\circle*{8}}
\put(120,20){\circle*{8}}
\put(160,20){\circle*{8}}
\put(200,20){\circle*{8}}
\put(  0,0){\makebox(0,0)[b]{$1$}}
\put( 40,0){\makebox(0,0)[b]{$2$}}
\put( 80,0){\makebox(0,0)[b]{$3$}}
\put(120,0){\makebox(0,0)[b]{$4$}}
\put(160,0){\makebox(0,0)[b]{$5$}}
\put(200,0){\makebox(0,0)[b]{$6$}}
\qbezier(  0,20)( 60, 80)(120,20)
\qbezier( 40,20)(100, 80)(160,20)
\qbezier( 80,20)(140, 80)(200,20)
\end{picture}
\caption{A perfect matching}\label{figure:matching}
\end{center}
\end{figure}

A permutation $(\sigma_{1},\sigma_{2},\hdots,\sigma_{n-1},\sigma_{n})$ 
which arises from a partition of $[n]$ into $2$-element blocks is called a \defterm{perfect matching} or a \defterm{$1$-factor}.
We say that the points $\sigma_{2i-1}$ and $\sigma_{2i}$ are \defterm{connected} to each other in this perfect matching $\sigma$.
We can express a perfect matching graphically by
arranging the lattice points $1$, $\dots$, $n$ along the $x$-axis in the plane 
and representing the edges $(\sigma_{2i-1},\sigma_{2i})$ by curves in the upper half plane.
Two edges $(\sigma_{2i-1},\sigma_{2i})$ and $(\sigma_{2j-1},\sigma_{2j})$ in $\sigma$ will be said to be \defterm{crossed}
if the corresponding edges intersect in such an embedding.
It is known that
$\sgn\sigma$ agrees with $(-1)^k$
where $k$ denotes the number of crossed pairs of edges in $\sigma$.
We write ${\cal F}_n$ for the set of perfect matchings of $[n]$.
For an example, the graphical representation of the perfect matching $\sigma=\{(1,4),(2,5),(3,6)\}\in{\cal F}_6$
is Figure~\ref{figure:matching} bellow,
and its sign is $-1$.

For each $\pi\in\frak{S}_n$, put $A^{\pi}=(a_{\pi(i)\pi(j)})$.
From the  definition above it is easy to see that 
\begin{equation}
\label{permutation}
\Pf(A^{\pi})=\sgn\pi\,\Pf(A).
\end{equation}
It is a well-known fact that the following identities hold.
For any skew symmetric $2n$ by $2n$ matrix $A$ and any $2n$ by $2n$ matrix $B$
we have
\begin{align}
&\Pf(A)^2=\det(A),\label{eq:square}\\
&\Pf\left(BA{}^tB\right)=\det(B)\Pf(A).
\nonumber
\end{align}
The first identity is fundamental
and we may use it implicitly hereafter.
The reader can prove it by the exterior algebra,
or can find a combinatorial proof in \cite{Ste}.
The second identity is a special case of Theorem~\ref{msf}.

Let $a_{ij}$ be a fixed element of a given square matrix $A$,
and denote by $(A;i,j)$ the square sub-matrix obtained by removing the $i$th row and $j$th column of $A$.
That is to say, we can write $(A;i,j)=A^{\overline{\{i\}}}_{\overline{\{j\}}}$ in the above notation,
where $\overline I$ stands for the complementary set of $I$.
The determinant of $(A;i,j)$ is called a \defterm{minor} corresponding to $a_{ij}$,
and the number $(-1)^{i+j}\det(A;j,i)$ is called a $(i,j)$-cofactor of $A$.
Here the terminology ``minor'' is used in a wider sense.
The  cofactor matrix $\widetilde A$ of $A$ is 
the matrix whose $(i,j)$-entry is the $(i,j)$-cofactor of A.
Then the following theorem is due to Jacobi.
\begin{thm}
\label{Lewis_Carroll}
Let $A$ be an $n$ by $n$ matrix
and $\widetilde A$ be its cofactor matrix.
Let $r\le n$ and $I,J\subseteq[n]$, $\sharp I=\sharp J=r$.
Then
\begin{equation}
\label{Lewis_Carroll_formula}
\det\widetilde A^{I}_{J}=(-1)^{|I|+|J|}(\det A)^{r-1}\;\det A^{\overline J}_{\overline I},
\end{equation}
where $\overline I,\overline J\subseteq[n]$ stand for the complements of $I,J$, respectively in $\subseteq[n]$.
Here we denote $|I|=\sum_{i\in I}i$.
\end{thm}
{\it Proof.}
Let $\Delta(i,j)=(-1)^{i+j}\,\det(A;j,i)$ denote the $(i,j)$-cofactor of $A$.
Then, by definition, the matrix of cofactors is
\begin{equation*}
\widetilde A
=\begin{pmatrix}
\Delta(i,j)
\end{pmatrix}
=\begin{pmatrix}
(-1)^{i+j}\det (A;j,i)
\end{pmatrix}.
\end{equation*}
Let $\overline I=\{p_1,\cdots,p_{n-r}\}$, $\overline J=\{q_1,\cdots,q_{n-r}\}$,
and let $\sigma : \bar{I} \to \bar{J}$ denote the order-preserving bijection 
which maps $p_k$ to $q_k$ for $k=1\dots n-r$.
Set $M=(M_{ij})$ to be the matrix defined by
\begin{equation*}
M_{ij}=\begin{cases}
\Delta(i,j) &\text{ if $i\in I$,}\\
\delta_{\sigma(i),j}&\text{ if $i\in\overline I$.}\\
\end{cases}
\end{equation*}
Then it is a direct simple algebra to see that the $(i,j)$-entry of the matrix 
$MA=B=(b_{ij})$ is
\begin{equation*}
b_{ij}=\begin{cases}
\delta_{ij}\,\det A&\text{ if $i\in I$,}\\
a_{\sigma(i),j}&\text{ if $i\in\overline I$.}\\
\end{cases}
\end{equation*}
Accordingly we have $\det B=(\det A)^{r}\det A^{\overline J}_{\overline I}$.
Meanwhile, it is easy to see that
\begin{equation*}
\det M
=(-1)^{\sum_{i\in\overline I}(i+\sigma(i))}\det \widetilde A^{I}_{J}
=(-1)^{|\overline I|+|\overline J|}\det \widetilde A^{I}_{J}
=(-1)^{|I|+|J|}\det \widetilde A^{I}_{J}
\end{equation*}
This proves the theorem.
$\Box$
\vskip .10in

\begin{ex}
\label{Dodgson}
We put $I=J=\{1,n\}\subset[n]$ in the formula above 
and obtain the Desnanot-Jacobi adjoint matrix theorem:
\begin{equation*}
\det M \det M^{2,\dots,n-1}_{2,\dots,n-1}
=\det M^{1,\dots,n-1}_{1,\dots,n-1}\det M^{2,\dots,n}_{2,\dots,n}
-\det M^{1,\dots,n-1}_{2,\dots,n} \det M^{2,\dots,n}_{1,\dots,n-1},
\end{equation*}
which is also called Dodgson's formula (or the Lewis Carroll formula).
For the details and the interesting story of the relations with the alternating sign matrices,
see \cite{B1}.
\end{ex}
\vskip .10in


Let $n$ be an even integer,
and let $A$ be a skew symmetric matrix of size $n$. 
For $1\le i<j\le n$, let $(A;\{i,j\},\{i,j\})$ denote the $(n-2)$ by $(n-2)$ skew symmetric sub-matrix
obtained by removing both the $i$th and $j$th rows and both the $i$th and $j$th columns of $A$,
i.e. $(A;\{i,j\},\{i,j\})=A^{\overline{\{i,j\}}}_{\overline{\{i,j\}}}$.
Let us define $\gamma(i,j)$ by
\begin{equation}
\gamma(i,j)=(-1)^{i+j-1} \Pf(A;\{i,j\},\{i,j\})
\end{equation}
for $1\le i<j\le n$.
We define the values of $\gamma(i,j)$ for $1\leq j\leq i\leq n$ so that $\gamma(j,i)=-\gamma(i,j)$ always holds.
Then the following expansion formula of Pfaffians along any row (resp. column) 
holds:
\begin{prop}
\label{expand}
Let $n$ be an even integer and $A=(a_{ij})$ be an $n$ by $n$ skew symmetric matrix. 
For any $i$, $j$ we have 
\begin{equation}
\label{expansion_row}
\delta_{ij}\Pf(A)
= \sum_{k=1}^{n} a_{kj}\gamma(k,i),
\end{equation}
\begin{equation}
\label{expansion_column}
\delta_{ij}\Pf(A)
= \sum_{k=1}^{n} a_{ik}\gamma(j,k).
\end{equation}
\end{prop}

Since $a_{ij}$ and $\gamma(i,j)$ are skew symmetric,
the reader sees immediately that the identities \thetag{\ref{expansion_row}} and \thetag{\ref{expansion_column}} are equivalent. Moreover, to prove the general case it is sufficient to show the case where $i=j=1$ in view of the formula 
 \thetag{\ref{permutation}}.
This case can be proved combinatorially from the definition 
 \thetag{\ref{def_pfaffian}} of Pfaffian.
If we multiply the both sides of \thetag{\ref{expansion_row}} by $\Pf(A)$ and use \thetag{\ref{eq:square}},
then we obtain
\begin{equation*}
\sum_{k=1}^{n} a_{ki}\gamma(k,j)\,\Pf(A)
=\delta_{ij}\,\left[\Pf(A)\right]^2=\delta_{ij}\,\det A.
\end{equation*}
Comparing this identity with the ordinary expansion of $\det A$,
we obtain  the following relation between
$\Delta(i,j)$ and $\gamma(i,j)$:
\begin{equation}
\label{Delta-gamma}
\Delta(i,j)=\gamma(j,i)\,\Pf(A).
\end{equation}
%
%
\begin{defi}
Let $n$ be an even integer.
Given a skew symmetric matrix $A$ of size $n$,
let us call $\gamma(i,j)$ a \defterm{copfaffian} corresponding to $a_{ij}$
(or \defterm{$(i,j)$-copfaffian}),
and let $\widehat A$ denote the skew symmetric matrix whose $(i,j)$-entry is $\gamma(i,j)$,
which we call the \defterm{copfaffian matrix} of $A$.
Note that \thetag{\ref{expansion_row}} and \thetag{\ref{expansion_column}} implies
\begin{equation}
{}^t\!\widehat A A=A\,{}^t\!\widehat A=\Pf(A)E_n,
\label{eq:copfaff}
\end{equation}
where $E_n$ denote the identity matrix of size $n$.
\end{defi}
\begin{ex}
Let $P_n(s,t)$ denotes the skew symmetric matrix,
whose $(i,j)$-entry is given by $s^{(i-1)\MOD 2+j\MOD 2}t^{j-i-1}$ for $1\leq i<j\leq n$,
where $x\MOD2$ stands for the remainder of $x$ divided by $2$.
In Lemma~7 of \cite{IW1},
we proved the formula
\begin{equation}
\label{product_formula}
\Pf(x_iy_j)_{1\leq i<j\leq n}=\prod_{i=1}^{[n/2]}x_{2i-1}\prod_{j=1}^{[n/2]}y_{2j}
\end{equation} for an even integer $n$.
From this formula, it is easy to see that the $(i,j)$-copfaffian of $P_n(s,t)$ is $(-1)^{j-i-1}s^{j-i-1}t^{(i-1)\MOD2+j\MOD2}$.
If $I=\{i_1,i_2,\hdots,i_{2r-1},i_{2r}\}_{<}$,
then the formula \thetag{\ref{product_formula}} also implies
\begin{equation}
\label{copfaff}
\Pf\left[P_n(s,t)^{I}_{I}\right]
=s^{\sum_{k=1}^{2r}(i_k-k)\MOD2}t^{\sum_{k=1}^{2r}(-1)^{k}i_{k}-r}.
\end{equation}
\end{ex}
The following result is considered as a Pfaffian version of Jacobi's formula.
\begin{thm}
\label{Lewis-Carroll_pfaffian}
Let $n$ be an even integer,
 and let $A$ be an $n$ by $n$ skew symmetric matrix.
Then, for any $I\subseteq[n]$ such that $\sharp\,I=2r$, we have
\begin{equation}
\label{LCpfaffian}
\Pf\big[(\widehat A)_{I}^{I}\big]=(-1)^{|I|-r}\left[\Pf (A)\right]^{r-1}\Pf(A_{\overline I}^{\overline I}).
\end{equation}
In particular, we have $\widehat{\widehat A}= (\Pf A)^{m-2} A$ with $n=2m$.
\end{thm}
{\it Proof.}
Let $\widetilde A=(\Delta(i,j))$ denote the matrix of the cofactors of $A$.
From \thetag{\ref{Delta-gamma}} we have $\widetilde A=\Pf(A)\,\widehat A$,
thus $(\widetilde A)^{I}_{I}=\Pf(A)\,(\widehat A)^{I}_{I}$.
It follows that
\begin{equation*}
\det(\widetilde A)^{I}_{I}
=\left[\Pf(A)\right]^{2r}\,\det(\widehat A)^{I}_{I}
=(\det A)^{r}\,\det(\widehat A)^{I}_{I}. 
\end{equation*}
On the other hand, Theorem~\ref{Lewis_Carroll} implies that
$\det(\widetilde A^{I}_{I})=(\det A)^{2r-1}\det A^{\overline I}_{\overline I}$. 
Comparing these two identities, we obtain
\begin{equation*}
\det(\widehat A)^{I}_{I}=(\det A)^{r-1} \det A^{\overline I}_{\overline I}. 
\end{equation*}
By taking the square root of both sides of this identity, we obtain 
\begin{equation}
\label{plusminus}
\Pf\big(\widehat A^{I}_{I}\big)
=\pm\left[\Pf(A)\right]^{r-1}\Pf\big(A^{\overline I}_{\overline I}\big). 
\end{equation}
Next we need to show that the signature in \thetag{\ref{plusminus}} does not depend of $A$.
Since the both sides of \thetag{\ref{plusminus}} are polynomials of the entries of $A$,
their ratio $\left[\Pf(A)\right]^{r-1}\Pf\big(A^{\overline I}_{\overline I}\big)/\Pf\big(\widehat A^{I}_{I}\big)$ is a rational function of them.
But this rational function can take only two values $\pm1$,
it must be a constant, i.e. independent of the entries of $A$.
To finish the proof we have to determine the sign. 
We substitute
\begin{equation}
S_{n}=P_{n}(1,1)=\begin{pmatrix}
 0&1&\dots&1\\
-1&0&\dots&1\\
\vdots&\vdots&\ddots&\vdots\\
-1&-1&\dots&0
\end{pmatrix}
\end{equation}
in the both sides of \thetag{\ref{plusminus}}.
Since $\widehat {S_n}=(\widehat S_{ij})$
with $\widehat S_{ij}=(-1)^{i+j-1}$ for $i<j$,
Applying \thetag{\ref{copfaff}}, we obtain $\Pf\widehat S_I=(-1)^{|I|-r}$ and $\Pf S_{\overline I}=1$,
which gives the desired sign.
Let $f(A)$ and $g(A)$ denote the left and right hand side of \thetag{\ref{LCpfaffian}}.
Then we have $(f+g)(f-g) = 0$ in the polynomial ring in the entries of $A$.
Since $f(S_n) + g(S_n) \neq 0$, 
we conclude that $f = g$.
This proves the theorem.
$\Box$

Given a skew symmetric matrix $A$,
we write $A(i_1,i_2,\dots,i_{2k})$ for
 $A_{i_1,i_2,\dots,i_{2k}}^{i_1,i_2,\dots,i_{2k}}$. 
\begin{ex}
\label{Pfaffian_Dodgson}
Given a skew symmetric matrix $A$ of size $n$,
take $I=\{1,2,3,4\}$ in the theorem,
then we obtain a formula which reads
\begin{align*}
\Pf(A)\Pf\left(A(5,\dots,n)\right)
&=\Pf\left(A({3,4,5,\dots,n})\right)\Pf\left(A({1,2,5,\dots,n})\right)\\
&\qquad-\Pf\left(A({2,4,5,\dots,n})\right)\Pf\left(A({1,3,5,\dots,n})\right)\\
&\qquad+\Pf\left(A({2,3,5,\dots,n})\right)\Pf\left(A({1,4,5,\dots,n})\right).
\end{align*}
This may be regarded as a Pfaffian version of Dodgson's identity given in Example~\ref{Dodgson}.
\end{ex}

We give some examples of the copfaffian matrices.
Let $n=2r$.
If $S_{n}=(S_{ij})$ with $S_{ij}=1$ for $i<j$,
then we have $\widehat {S_n}=(\widehat S_{ij})$ with $\widehat S_{ij}=(-1)^{i+j-1}$ for $i<j$ as obtained in the above proof.
Put $T_{n}=\left(T_{ij}\right)=P_{n}(0,1)$ and $\widehat{T_{n}}=\left({\widehat T}_{ij}\right)$,
then
\begin{align*}
&T_{ij}=\begin{cases}
1&\text{ if $1\leq i<j$, and $i$ and $j-i$ are both odd,}\\
0&\text{otherwise.}
\end{cases}\\
&{\widehat T}_{ij}=\begin{cases}
1&\text{ if $j=i+1$,}\\
0&\text{otherwise}.
\end{cases}
\end{align*}
For instance, we have 
\begin{equation*}
T_{4}=\begin{pmatrix}
0&1&0&1\\
-1&0&0&0\\
0&0&0&1\\
-1&0&-1&0\\
\end{pmatrix},
\qquad
{\widehat T}_{4}=\begin{pmatrix}
0&1&0&0\\
-1&0&1&0\\
0&-1&0&1\\
0&0&-1&0\\
\end{pmatrix}.
\end{equation*}
Let $E_{m}$ denote the identity matrix of size $m$,
and let $O_{m,n}$ denote the $m$ by $n$ zero matrix.
When $m=n$, we simply write $O_{m}$ for $O_{m,m}$.
Let $J_{m}$ denotes the symmetric matrix of size $m$ defined by
\begin{equation*}
J_{m}=\begin{pmatrix}
     0&\hdots&     0&     1\\
     0&\hdots&     1&     0\\
\vdots&\rdots&\vdots&\vdots\\
     1&\hdots&     0&     0\\
\end{pmatrix}.
\end{equation*}
We let $n=2m$ and put $K_{n}=\begin{pmatrix} O_m&J_{m}\\ -J_{m}&O_{m}\end{pmatrix}$ and $L_{n}=\begin{pmatrix} O_{m}&E_{m}\\ -E_{m}&O_{m}\end{pmatrix}$.
Then it is easy to see that ${}^{t}\kern-1pt K_{n}=-K_{n}$, ${}^{t}\kern-1pt L_{n}=-L_{n}$, $K_{n}^{2}=L_{n}^{2}=-I_{n}$, $\Pf(K_{n})=1$ and $\Pf(L_{n})=(-1)^{\frac{m(m-1)}2}$.
From Cramer's formula and \thetag{\ref{Delta-gamma}},
we have $\widehat A=\Pf(A)\,{}^{t}\hskip-2pt A^{-1}$ for a non-singular matrix $A$,
which immediately implies $\widehat{K_{n}}=K_{n}$ and $\widehat{L_{n}}=(-1)^{m(m-1)/2}L_{n}$.

We next state a Pfaffian analogue of the Pl\"ucker relations (or known as the Grassmann-Pl\"ucker relations for determinants) and make a remark on a  
relation with the Lewis Carroll formula.
It is an algebraic identity of degree two describing the relations among several subpfaffians. We point out here that 
this identity has been proved in \cite{Hi} and in \cite{DW} in the framework of an exterior algebra.
\begin{thm}
\label{Plucker}
Suppose $m$, $n$ are odd integers.
Let $A$ be an $(m+n)\times(m+n)$ skew symmetric matrix.
Fix a sequence of integers $I=\{i_1,i_2,\dots,i_m\}_{<}\subseteq[m+n]$ such that $\sharp\,I=m$.
Denote the complement of $I$ by $\overline I = \{k_1,k_2,\dots,k_n\}_{<}\subseteq[m+n]$ which has the cardinality $n$. 
Then the following relation holds.
\begin{equation}
\sum_{j=1}^{m}(-1)^{j-1}\Pf\left(A_{I\setminus\{i_j\}}^{I\setminus\{i_j\}}\right)\Pf\left(A_{\{i_j\}\cup\overline I}^{\{i_j\}\cup\overline I}\right)
=\sum_{j=1}^{n}(-1)^{j-1}\Pf\left(A_{I\cup\{k_j\}}^{I\cup\{k_j\}}\right)\Pf\left(A_{\overline I\setminus \{k_j\}}^{\overline I\setminus \{k_j\}}\right).
\end{equation}
\end{thm}
{\it Proof.}
We only use the expansion formula of a Pfaffian given in Proposition~\ref{expand}.
In fact, if we expand $\Pf(A_{\{i_j\}\cup\overline I}^{\{i_j\}\cup\overline I})$ along the  $i_j$th row/column on the left-hand side
and also expand $\Pf(A_{I\cup\{k_j\}}^{I\cup\{k_j\}})$ along the $k_j$th row/column on the right-hand side,
and compare with each other, then it is immediate to see the desired equality.
This identity is also proved directly from the definition \thetag{\ref{def_pfaffian}} of a Pfaffian
by using the notion of matching and related combinatorics.
$\Box$
\vskip .10in

The formula in the following assertion, which is called by the basic 
identity in \cite{Kn}, is regarded as a special case of the Pl\"ucker relations
above.
\begin{cor}
\label{basic_identity}
Let $A$ be a skew symmetric matrix of size $N$.
Let $I=\{i_1,i_2,\dots,i_{2k}\}$ be a subset of $[N]$.
Take an integer $l$ which satisfies $2k+2l\leq N$.
Then we have 
\begin{align}
&\Pf(A(1,2,\dots,2l))
\Pf(A(i_1,i_2,\dots,i_{2k},1,\dots,2l))\nonumber\\
=& \sum_{j=2}^{2k}(-1)^{j}
\Pf(A(1,2,\dots,2l,i_1,i_{j}))
\Pf(A(i_2,\dots,\widehat i_{j},\dots,i_{2k},1,\dots,2l)).
\end{align}
\end{cor}
{\it Proof.}
Given a skew symmetric matrix $A=(a_{ij})_{1\leq i,j\leq N}$ of size $N$ and a subset $I=\{i_1,i_2,\dots,i_{2k}\}$,
we consider the skew symmetric matrix $B=(b_{ij})_{1\leq i,j\leq 2k+4l}$ of size $2k+4l$,
whose $(i,j)$-entry $b_{ij}$ is equal to $a_{p_{i}p_{j}}$,
where the sequence $\{p_{i}\}_{i=1}^{2k+4l}$ is determined by
\begin{align*}
\begin{cases}
p_{\nu}=\nu&\text{ if $1\leq\nu\leq2l$,}\\
p_{2l+\nu}=i_{\nu}&\text{ if $1\leq\nu\leq2k$,}\\
p_{2k+2l+\nu}=\nu&\text{ if $1\leq\nu\leq2l$.}
\end{cases}
\end{align*}
Now, apply Theorem~\ref{Plucker} to $B$ with $m=2l+1$, $n=2k+2l-1$, $I=\{1,2,\dots,2l+1\}$ and $\overline I=\{2l+2,2l+3,\dots,2k+4l\}$.
Then, since  each summand on the left-hand side disappears except for the case
$j=2l-1\,(i_j=i_1)$,
the desired identity immediately follows from the identity in Theorem~\ref{Plucker}.
$\Box$
\vskip .10in

\begin{rem}
If we put $k=2$ in this corollary,
then the identity is nothing but the identity in Example~\ref{Pfaffian_Dodgson}.
This implies the basic identity partially covers the Lewis Carroll formula
for Pfaffians.
\end{rem}

%% file: lmsf3.tex
%
%
%
%
%
%

\section{Summation formulas of Pfaffians}

In this section we review the summation formulas of Pfaffians.
We restate the theorems in \cite{IW1} 
and,
in the next section, we will give 
certain combinatorial proofs of Theorem~\ref{msf}, \ref{msf2} and \ref{msf3}.
In \cite{IW1} we gave algebraic proofs of these theorems,
whereas,
in this paper we use the lattice paths combined with the Pfaffian version of Jacobi's formula (i.e. Theorem~\ref{Lewis-Carroll_pfaffian})
to prove these theorems.
Especially we will find that the Pfaffian version of Jacobi's formula is useful to simplify 
our lattice path proof,
comparing with the combinatorial discussion given in \cite{Ste},
and give more insights to explain these formulas.
Our proofs of the theorems described here will be postponed until the 
next section.
%
%
%
%
\begin{lem}
\label{pf_lem3}
Let $n$, $m$ and $M$ be nonnegative integers,
and let $N=2N'$ be an even integer.
Let $A$ (resp. $B$) be a skew symmetric matrix of size $M$ (resp. $N$) such that $B$ is non-singular.
Let $T_{11}$, $T_{12}$, $T_{21}$ and $T_{22}$ be an $m$ by $n$, $m$ by $N$, $M$ by $n$ and $M$ by $N$ rectangular matrix, respectively.
Then
\begin{align}
&\Pf(B)^{-1}
\Pf\begin{pmatrix}
 T_{12}B{}^t\! T_{12} & T_{12}B{}^t\! T_{22} & T_{11}J_n\\
 T_{22}B{}^t\! T_{12} &A+T_{22}B\,{}^t\! T_{22}&T_{21} J_{n}\\
-J_n {}^t\! T_{11}&-J_{n}{}^t\! T_{21}&O_{n}
\end{pmatrix}
\nonumber\\
&=\Pf\begin{pmatrix} 
   O_m   &   O_{m,M} & T_{12} J_N & T_{11} J_n\\
O_{M,m}  & A       &   T_{22} J_N & T_{21} J_n\\
- J_N {}^t\! T_{12} & - J_N {}^t\! T_{22}& \frac1{\Pf(B)}J_N {}^t\!\widehat{B} J_N& O_{N,n}\\
- J_n {}^t\! T_{11} & - J_n {}^t\! T_{21}& O_{n,N} & O_{n}
\end{pmatrix}
\nonumber\\
&=\Pf\begin{pmatrix} 
   J_M {}^t\!A J_M   &   O_{M,m} & J_M T_{21}  & J_M T_{22}\\
O_{m,M}        &   O_m     & J_m  T_{11} & J_m T_{12}\\
-{}^t\! T_{21} J_M & -{}^t\! T_{11} J_m& O_{n} & O_{n,N}\\
-{}^t\! T_{22} J_M & -{}^t\! T_{12} J_m& O_{N,n} & \frac1{\Pf(B)}\widehat{B}
\end{pmatrix}
\end{align}
\end{lem}
{\it Proof.}
The first identity follows from the following matrix identity
\begin{align*}
&\begin{pmatrix} 
   O_m   &   O_{m,M} & T_{12} J_N & T_{11} J_n\\
O_{m,n}  & A       &   T_{22} J_N & T_{21} J_n\\
- J_N\, {}^t\! T_{12} & - J_N\, {}^t\! T_{22}& \frac1{\Pf(B)}J_N\, {}^t\!\widehat{B} J_N& O_{N,n}\\
- J_n\, {}^t\! T_{11} & - J_n\, {}^t\! T_{21}& O_{n,N} & O_{n}
\end{pmatrix}
\begin{pmatrix}
E_m&O_{m,M}&O_{m,N}&O_{m,n}\\
O_{M,m}& E_M & O_{M,N}& O_{M,n} \\
J_{N}B\,{}^t\!T_{12} & J_{N}B\,{}^t\!T_{22} & E_{N} & O_{N,n} \\
O_{n,m} & O_{n,M} & O_{n,N} & E_{n}\\
\end{pmatrix}\\
&=\begin{pmatrix} 
T_{12}B\,{}^t\!T_{12}   &   T_{12}B\,{}^t\!T_{22} & T_{12} J_N & T_{11} J_n\\
T_{22}B\,{}^t\!T_{12}  & A+ T_{22}B\,{}^t\!T_{22}      &   T_{22} J_N & T_{21} J_n\\
O_{M,m} & O_{N,n}& \frac1{\Pf(B)}J_N\, {}^t\!\widehat{B} J_N& O_{N,n}\\
- J_n\, {}^t\! T_{11} & - J_n\, {}^t\! T_{21}& O_{n,N} & O_{n}
\end{pmatrix},
\end{align*}
which follows from \thetag{\ref{eq:copfaff}}.
By taking the determinants of both sides of this equation
and using $\Pf(J_N\, {}^t\!\widehat{B} J_N)=\Pf(\widehat B)=\Pf(B)^{-1}$ from \thetag{\ref{LCpfaffian}},
one obtains
\begin{equation*}
\Pf\begin{pmatrix}
 T_{12}B{}^t\! T_{12} & T_{12}B{}^t\! T_{22} & T_{11}J_n\\
 T_{22}B{}^t\! T_{12} &A+T_{22}B\,{}^t\! T_{22}&T_{21} J_{n}\\
-J_n {}^t\! T_{11}&-J_{n}{}^t\! T_{21}&O_{n}
\end{pmatrix}
=\pm\Pf(B)
\Pf\begin{pmatrix} 
   O_m   &   O_{m,M} & T_{12} J_N & T_{11} J_n\\
O_{M,m}  & A       &   T_{22} J_N & T_{21} J_n\\
- J_N {}^t\! T_{12} & - J_N {}^t\! T_{22}& \frac1{\Pf(B)}J_N {}^t\!\widehat{B} J_N& O_{N,n}\\
- J_n {}^t\! T_{11} & - J_n {}^t\! T_{21}& O_{n,N} & O_{n}
\end{pmatrix}
\end{equation*}
Let $f(A,B,T_{11},T_{12},T_{21},T_{22})$ (resp. $g(A,B,T_{11},T_{12},T_{21},T_{22})$)
denote the left-hand side (resp. the right-hand side) of this equation.
Then we have $(f-g)(f+g)=0$ in the polynomial ring of the entries of $A$, $B$, $T_{11}$, $T_{12}$, $T_{21}$, $T_{22}$.
Comparing the coefficients of the entries of $A$, we conclude that $f+g\neq0$,
which implies $f=g$.
This shows that the signature does not depend on $A$, $B$ and $T_{ij}$.
The other identity is proved similarly.
$\Box$
\bigbreak

We restate Theorem~1 of \cite{IW1} in the following form:
%
%
\begin{thm}
\label{msf}
Let $m$ and $N=2N'$ be even integers such that $m\le N$.
Let $T=(t_{ik})_{1\le i\le m, 1\le k\le N}$ be an $m$ by $N$ rectangular matrix.
Let $A=(a_{ij})_{1\le i,j\le N}$ be a non-singular skew-symmetric matrix of size $N$
and let $\widehat A$ denote its copfaffian matrix.
Then
\begin{align}
\label{eq_msf}
&\sum_{I\in\binom{[N]}{m}}
\Pf(A^I_I) \det(T_I)
=\Pf(Q)\nonumber\\
&=\Pf(A)\Pf\begin{pmatrix}
   O_m   &   T J_N \\
-J_N {}^tT & \frac1{\Pf(A)}J_N{}^t\!\widehat{A}J_N \\
\end{pmatrix}
=\Pf(A)\Pf\begin{pmatrix}
   O_m  & J_m T \\
-{}^tT J_m    & \frac1{\Pf(A)}\widehat{A} \\
\end{pmatrix}.
\end{align}
Here $Q=(Q_{ij})=TA\,{}^t\! T$, and its entries are given by
\begin{equation}
\label{pf_msf}
Q_{ij}=\sum_{1\le k<l\le N} a_{kl} \det(T^{ij}_{kl}),
\qquad(1\leq i,j\leq m).
\end{equation}
\end{thm}
The first identity, i.e. $\sum\Pf(A^I_I) \det(T_I)=\Pf(Q)$,
holds even if $N$ is not even.
When $m$ is odd,
we can immediately derive a similar formula from the case where $m$ is even.
So we only treat even cases in this paper.
If we take $m=N=2r$ and $A=K_{m}$ in \thetag{\ref{eq_msf}},
then $\det(T)=\Pf(K_{m})\det(T)=\Pf(T\,K_{m}{}^t\!T)$.
This means that every determinant of even degree can be represented by a Pfaffian of the {\it same} degree.

When $m$ and $N$ are even integers such that $0\leq m\leq N$,
and $X$ and $Y$ are $m$ by $N$ rectangular matrices,
taking $A={}^{t}\!X K_{m} X$ and $T=Y$ in Theorem~{\ref{msf}},
 we obtain the following corollary, which is the so-called Cauchy-Binet formula.
For another proof which use Theorem~\ref{msf} also, see \cite{IOW}.  
%
%
\begin{cor}
\label{cauchy-binet}
Assume $m\le N$,
and let $X=(x_{ij})_{1 \le i \le m, 1 \le j \le N}$ and $Y=(y_{ij})_{1 \le i \le m, 1 \le j \le N}$ be any $m$ by $N$ matrices.
Let $A=(a_{ij})_{1 \le i \le N, 1 \le j \le N}$ be any $N$ by $N$ matrix.
Then
\begin{equation}
\label{eq_cauchy-binet}
\sum_{K\in\binom{[N]}{m}}
 \det(X_K)
 \det(Y_K)
= \det \left( X \trans Y \right).
\end{equation}
Especially, using \thetag{\ref{eq_cauchy-binet}} twice,
we obtain the following general version:
\begin{equation}
\label{eq_gene_cauchy-binet}
\sum_{I,J\in\binom{[N]}{m}}
 \det(A^I_J)
\det(X_I)
 \det(Y_J)
= \det \left( X A\, \trans Y \right).
\end{equation}
\end{cor}

The following theorem gives a minor summation formula,
in which the index set $I$ of a minor in the sum
always includes some fixed column index set, say $\{1,2,\hdots,n\}$.
(See \cite{IW1} and \cite{Ste}.) 
%
%
\begin{thm}
\label{msf2}
Let $m$, $n$ and $N$ be positive integers such that $m-n$ and $N$ are even and $0\le m-n\le N$.
Let $A=(a_{ij})_{1\le i,j\le N}$ be a nonsingular skew-symmetric matrix of size $N$.
Let $T=(t_{ij})_{1\le i\le m, 1\le j\le n+N}$ be an $m$ by $(n+N)$ rectangular matrix.
Write the sets of column indices as 
$R^0=\{1,\dots,n\}$ and $R=\{n+1,\dots,n+N\}$.
Then
\begin{align}
\label{minor_sum}
&\sum_{I\in\binom{R}{m-n}}
\Pf(A^I_I) \det(T_{R^0\uplus I})=\Pf\begin{pmatrix} Q&T_{R^0} J_n\\ -J_n{}^t\!T_{R^0}&O_n\end{pmatrix}
\nonumber\\
&=\Pf(A)\Pf\begin{pmatrix}
   O_m   &   T_{R} J_N & T_{R^0} J_n \\
- J_N {}^t\!T_{R} & \frac1{\Pf(A)}J_N {}^t\!\widehat{A} J_N& O_{N,n}\\
- J_n {}^t\!T_{R^0} & O_{n,N} & O_{n}
\end{pmatrix},
\end{align}
where $Q$ is the $m$ by $m$ skew-symmetric matrix defined by $Q=T_{R}A\,{}^t\!T_{R}$, i.e.,
\begin{equation}
\label{Q}
Q_{ij}=\sum_{1\le k<l\le N} a_{kl} \det(T^{ij}_{kl}),
\qquad(1\le i,j\le m).
\end{equation}
\end{thm}
We shall later restate this theorem (and a proof) as Theorem~\ref{lmsf2}
in a combinatorial description using lattice paths quite naturally.


The following theorem shows
 a minor summation formula for both the rows and columns.
%
%
\begin{thm}
\label{msf3}
Let $M$ and $N$ be even integers such that $M\le N$.
Let $T=(t_{ij})_{1\le i\le M,1\le j\le N}$ be any $M$ by $N$ rectangular matrix,
and let $A=(a_{ij})_{1\le i,j\le M}$ (resp. $B=(b_{ij})_{1\le i,j\le N}$) be a nonsingular skew-symmetric matrix of size $M$ (resp. size $N$).
Then
\begin{align}
&\sum_{r=0}^{\lfloor M/2\rfloor}
z^{2r}\sum_{I\in\binom{[M]}{2r}}
\sum_{J\in\binom{[N]}{2r}}
\Pf(A^I_I)\Pf(B^J_J) \det(T^I_J)
=\Pf(A)\Pf\left[\frac1{\Pf(A)}\widehat A+z^2Q\right]
\nonumber\\
&
=\Pf\begin{bmatrix}
J_M\,{}^t\kern-2pt AJ_M&J_M\\
-J_M&z^2Q
\end{bmatrix}
=\Pf(A)\Pf(B)
\Pf\left[\begin{matrix}
\frac1{\Pf(A)}\widehat A&zTJ_N\\
-zJ_N{}^t\! T&\frac1{\Pf(B)}J_N\,{}^t\!\widehat B\,J_N
\end{matrix}\right]
\nonumber\\
&=\Pf(A)\Pf(B)
\Pf\left[\begin{matrix}
\frac1{\Pf(A)}J_n\,{}^t\!\widehat A J_n&zJ_n T\\
-z{}^t\! T J_n&\frac1{\Pf(B)}\widehat B
\end{matrix}\right]
\end{align}
where $Q=TB\,{}^tT$ and $\lfloor x\rfloor$ denotes the
greatest integer that does not exceed $x$.
\end{thm}
We also have the 
%
%
\begin{cor}
\label{msf4}
Let $M$ and $N$ be nonnegative integers such that $M\le N$.
Let $T=(t_{ij})$ be an $M$ by $N$ rectangular matrix.
Let $B=(b_{ij})_{0\leq i,j\leq N}$ be a skew-symmetric matrix of size $(N+1)$.
\begin{enumerate}
\item
If $M$ is odd and $A=(a_{ij})_{0\leq i,j\leq M}$ is a nonsingular skew-symmetric matrix of size $(M+1)$, then
\begin{align}
\sum_{r=0}^{\lfloor M/2\rfloor}z^{2r}
&\sum_{I\in\binom{[M]}{2r}}
\sum_{J\in\binom{[N]}{2r}}
\Pf(A_{I}^{I}) \Pf({B}_{J}^{J}) \det(T^{I}_{J})\nonumber\\
&+\sum_{r=0}^{\lfloor (M-1)/2\rfloor}z^{2r+1} 
\sum_{I\in\binom{[M]}{2r+1}}
\sum_{J\in\binom{[N]}{2r+1}}
\Pf(A_{\{0\}\cup I}^{\{0\}\cup I})
\Pf({B}_{\{0\}\cup J}^{\{0\}\cup J})\det(T^{I}_{J})\nonumber\\
&=\Pf(A)\Pf\left(\frac1{\Pf(A)}\widehat A+Q\right),
\end{align}
where $Q=(Q_{ij})_{0\leq i,j\leq M}$ is given by
\begin{equation}
Q_{ij}=\begin{cases}
0,&\qquad\text{ if }i=j=0,\\
z\sum_{1\le k\le N} b_{0k} t_{jk},&\qquad\text{ if }i=0\text{ and }1\le j\le M,\\
z\sum_{1\le k\le N} b_{k0} t_{jk},&\qquad\text{ if }j=0\text{ and }1\le i\le M,\\
z^2\sum_{1\le k<l\le N} b_{kl} \det(T^{ij}_{kl}),&\qquad\text{ if }1\le i,j\le M.\end{cases}
\end{equation}
\item
If $M$ is even and $A=(a_{ij})_{0\leq i,j\leq M+1}$ is a nonsingular skew-symmetric matrix of size $(M+2)$, then
\begin{align}
\sum_{r=0}^{\lfloor M/2\rfloor}
z^{2r} &\sum_{I\in\binom{[M]}{2r}}
\sum_{J\in\binom{[N]}{2r}}
\Pf(A_{I}^{I}) \Pf({B}_{J}^{J}) \det(T^{I}_{J})\nonumber\\
&+\sum_{r=0}^{\lfloor (M-1)/2\rfloor}z^{2r+1} 
\sum_{I\in\binom{[M]}{2r+1}}
\sum_{J\in\binom{[N]}{2r+1}}
\Pf(A_{\{0\}\cup I}^{\{0\}\cup I})
\Pf({B}_{\{0\}\cup J}^{\{0\}\cup J})\det(T^{I}_{J})\nonumber\\
&=\Pf(A)\Pf\left(\frac1{\Pf(A)}\widehat A+Q\right),
\end{align}
where $Q=(Q_{ij})_{0\leq i,j\leq M+1}$ is given by
\begin{equation}
Q_{ij}=\begin{cases}
0&\qquad\text{ if $i=j=0$},\\
z\sum_{1\le k\le N} b_{0k} t_{jk},&\qquad\text{ if $i=0$ and $1\leq j\leq M$,}\\
z\sum_{1\le k\le N} b_{k0} t_{jk},&\qquad\text{ if $j=0$ and $1\leq i\leq M$,}\\
z^2\sum_{1\le k<l\le N} b_{kl} \det(T^{ij}_{kl}),&\qquad\text{ if $1\leq i,j\leq M$,}\\
0&\qquad\text{ if $i=M+1$ or $j=M+1$.}\\
\end{cases}
\end{equation}
\end{enumerate}
\end{cor}

%
%
\begin{thm}
\label{msf4}
Let $m$, $n$, $M$ and $N$ be non-negative integers
 such that $M$, $N$ and $m-n$ are even.
We put $R^0=\{1,\dots,m\}$, $S^0=\{1,\dots,n\}$, $R=\{m+1,\dots,m+M\}$ and S=\{n+1,\dots,n+N\}.
Let $T=(t_{ij})_{1\le i\le m+M,1\le j\le n+N}$ be any $(m+M)$ by $(n+N)$ rectangular matrix,
and let $A=(a_{ij})_{1\le i,j\le M}$ (resp. $B=(b_{ij})_{1\le i,j\le N}$) be a nonsingular skew-symmetric matrix of size $M$ (resp. size $N$).
Then
\begin{align}
&\sum_{{\max(m,n)\leq r\leq\min(m+M,n+N)}\atop{\text{$r-\max(m,n)$ is even.}}}
z^{r}\sum_{I\in\binom{R}{r-m}}\sum_{J\in\binom{S}{r-n}}
\Pf(A^I_I)\Pf(B^J_J) \det(T^{R^0\uplus I}_{S^0\uplus J})
\nonumber\\
&=\Pf(A)\Pf(B)\Pf\begin{pmatrix} 
   O_m   &   O_{m,M} & zT^{R^0}_{S} J_N & zT^{R^0}_{S^0}J_n\\
O_{M,m}  & \frac1{\Pf(A)}\widehat{A}       &   zT^{R}_{S} J_N & zT^{R}_{S^0} J_n\\
-zJ_N{}^t\!T^{R^0}_{S}  & -zJ_N{}^t\!T^{R}_{S}& \frac1{\Pf(B)}J_N {}^t\!\widehat{B} J_N& O_{N,n}\\
-zJ_n{}^t\!T^{R^0}_{S^0}& -zJ_n{}^t\!T^{R}_{S^0}& O_{n,N} & O_{n}
\end{pmatrix}
\end{align}
\end{thm}
%
%
\begin{cor}
Let $m$, $n$, $M$ and $N$ be nonnegative integers such that $m-n$ is even and $M\leq N$.
We put $R^0=\{1,\dots,m\}$, $S^0=\{1,\dots,n\}$, $R=\{m+1,\dots,m+M\}$ and S=\{n+1,\dots,n+N\}. 
Let $T=(t_{ij})_{1\le i\le m+M,1\le j\le n+N}$ be any $(m+M)$ by $(n+N)$ rectangular matrix 
and let $B=(b_{ij})_{1\le i,j\le N+1}$ be any skew-symmetric matrix of size $(N+1)$.
\begin{enumerate}
\item
If $M$ is odd and $A=(a_{ij})_{1\leq i,j\leq M+1}$ is a nonsingular skew-symmetric matrix of size $(M+1)$, then
\begin{align}
&\sum_{{\max(m,n)\leq r\leq\min(m+M,n+N)}\atop{r-\max(m,n)\text{ is even.}}}
z^{r}\sum_{I\in\binom{R}{r-m}}\sum_{J\in\binom{S}{r-n}}
\Pf(A^I_I)\Pf(B^J_J) \det(T^{R^0\uplus I}_{S^0\uplus J})
\nonumber\\
&+\sum_{{\max(m,n)\leq r\leq\min(m+M,n+N)}\atop{r-\max(m,n)\text{ is odd.}}}
z^{r}\sum_{I\in\binom{R}{r-m}}\sum_{J\in\binom{S}{r-n}}
\Pf(A^{I\uplus\{M+1\}}_{I\uplus\{M+1\}})
\Pf(B^{J\uplus\{N+1\}}_{J\uplus\{N+1\}})
\det(T^{R^0\uplus I}_{S^0\uplus J})
\nonumber\\
&=\Pf(A)\Pf\begin{pmatrix} 
Q^{11}&Q^{12}& zT^{R^0}_{S^0}J_n\\
-{}^t\!Q^{12} & \frac1{\Pf(A)}\widehat{A}+Q^{22}& z {\overline T}^{R}_{S^0} J_n\\
-zJ_n{}^t\!T^{R^0}_{S^0}  &-zJ_n{}^t\!{\overline T}^{R}_{S^0}  & O_n
\end{pmatrix},
\end{align}
where $Q^{11}=(Q^{11}_{ij})_{1\leq i,j\leq m}$,
 $Q^{12}=(Q^{12}_{ij})_{1\leq i\leq m,\:1\leq j\leq M+1}$,
$Q^{22}=(Q^{22}_{ij})_{1\leq i,j\leq M+1}$ is given by
\begin{align*}
&Q^{11}_{ij}=z^2\sum_{1\le k<l\le N} b_{kl} \det(T^{ij}_{n+k,n+l})
\qquad\text{ for $1\leq i,j\leq m$,}\\
&Q^{12}_{ij}=\begin{cases}
z^2\sum_{1\le k<l\le N} b_{kl} \det(T^{i,j+m}_{n+k,n+l})&\text{ if $1\leq i\leq m$ and $1\leq j\leq M$,}\\
z\sum_{1\le k\le N} b_{k,N+1}T^{i}_{n+k}&\text{ if $1\leq i\leq m$ and $j=M+1$,}\\
\end{cases}\\
&Q^{22}_{ij}
=\begin{cases}
z^2\sum_{1\le k<l\le N} b_{kl} \det(T^{m+i,m+j}_{n+k,n+l})&\text{ if $1\le i,j\le M$,}\\
z\sum_{1\le k\le N} b_{k,N+1} T^{m+i}_{n+k}&\text{ if $1\le i\le M$ and $j=M+1$,}\\
z\sum_{1\le k\le N} b_{N+1,k} T^{m+j}_{n+k}&\text{ if $i=M+1$ and $1\le j\le M$,}\\
0&\text{ if $i=j=M+1$,}
\end{cases}
\end{align*}
and $\overline {T}^{R}_{S^0}$ is the $(M+1)$ by $n$ matrix in which its first $M$ rows are the same as $T^{R}_{S^0}$
and the entries of its bottom row are all zero.
\item
If $M$ is even and $A=(a_{ij})_{1\leq i,j\leq M+2}$ is a nonsingular skew-symmetric matrix of size $(M+2)$, then
\begin{align}
&\sum_{{\max(m,n)\leq r\leq\min(m+M,n+N)}\atop{r-\max(m,n)\text{ is even.}}}
z^{r}\sum_{I\in\binom{R}{r-m}}\sum_{J\in\binom{S}{r-n}}
\Pf(A^I_I)\Pf(B^J_J) \det(T^{R^0\uplus I}_{S^0\uplus J})
\nonumber\\
&+\sum_{{\max(m,n)\leq r\leq\min(m+M,m+N)}\atop{r-\max(m,n)\text{ is odd.}}}
z^{r}\sum_{I\in\binom{R}{r-m}}\sum_{J\in\binom{S}{r-n}}
\Pf(A^{I\uplus\{M+2\}}_{I\uplus\{M+2\}})\Pf(B^{J\uplus\{N+1\}}_{J\uplus\{N+1\}}) \det(T^{R^0\uplus I}_{S^0\uplus J})
\nonumber\\
&=\Pf(A)\Pf\begin{pmatrix} 
Q^{11}&Q^{12}& zT^{R^0}_{S^0}J_n\\
-{}^t\!Q^{12} & \frac1{\Pf(A)}\widehat{A}+Q^{22}& z{T^*}^{R}_{S^0} J_n\\
-zJ_n{}^t\!T^{R^0}_{S^0}  &-zJ_n{}^t\!{T^*}^{R}_{S^0}  & O_n
\end{pmatrix},
\end{align}
where $Q^{11}=(Q^{11}_{ij})_{1\leq i,j\leq m}$,
$Q^{12}=(Q^{12}_{ij})_{1\leq i\leq m,\:1\leq j\leq M+2}$,
$Q^{22}=(Q^{22}_{ij})_{1\leq i,j\leq M+2}$ is given by
\begin{align*}
&Q^{11}_{ij}=z^2\sum_{1\le k<l\le N} b_{kl} \det(T^{ij}_{n+k,n+l})
\qquad\text{ for $1\leq i,j\leq m$,}\\
&Q^{12}_{ij}=\begin{cases}
z^2\sum_{1\le k<l\le N} b_{kl} \det(T^{i,j+m}_{n+k,n+l})
&\text{ if $1\leq i\leq m$ and $1\leq j\leq M$,}\\
0&\text{ if $1\leq i\leq m$ and $j=M+1$,}\\
z\sum_{1\le k\le N} b_{k,N+1}T^{i}_{n+k}&\text{ if $1\leq i\leq m$ and $j=M+2$,}\\
\end{cases}\\
&Q^{22}_{ij}
=\begin{cases}
z^2\sum_{1\le k<l\le N} b_{kl} \det(T^{m+i,m+j}_{n+k,n+l})&\text{ if $1\le i,j\le M$,}\\
z\sum_{1\le k\le N} b_{k,N+1} T^{m+i}_{n+k}&\text{ if $1\le i\le M$ and $j=M+2$,}\\
z\sum_{1\le k\le N} b_{N+1,k} T^{m+j}_{n+k}&\text{ if $i=M+2$ and $1\le j\le M$,}\\
0&\text{ otherwise,}
\end{cases}
\end{align*}
and ${T^*}^{R}_{S^0}$ is the $(M+2)$ by $n$ matrix in which its first $M$ rows are the same as $T^{R}_{S^0}$
and the entries of the last two rows are all zero.
\end{enumerate}
\end{cor}

%% file: lmsf4.tex
%
%
%
%
%
%

\section{Proofs by Lattice Paths}

In this section we give combinatorial proofs of the summation formulas of Pfaffians,
i.e.,
Theorem~\ref{msf}, \ref{msf2}, \ref{msf3},
which are stated in Section 3.
In \cite{O1} Okada gave this type of the formula related to a certain plane partition enumeration problem,
but his formula was a very special case, that is,  $A=S_{N}$, of ours.
In \cite{Ste} J.Stembridge gave a lattice path interpretation of this special summation formula,
and gave proofs from this point of view.
We follow his line in part and actually give a lattice path interpretation of 
our formulas and proofs from this viewpoint.
However, it is important to notice here that the Pfaffian analogue of the 
Lewis Carroll formula (Theorem~\ref{Lewis-Carroll_pfaffian}) makes 
possible the story clear. 
Thus, in this section, we provide an improved and a 
much simplified version of Stembridge's proof.
We may say our proofs are closer to Gessel-Viennot's original proofs in \cite{GV} than those given in \cite{Ste}.
We note that Stembrige's proof \cite{Ste} 
can be also generalized almost parallelly to proving
 Theorem~\ref{msf}, \ref{msf2}, \ref{msf3},
but we do not develop the proofs in this direction.

Now we review the basic terminology of lattice paths and fix notation.
We follow the basic terminology in \cite{GV} and \cite{Ste}.
Let $D=(V,E)$ be an acyclic digraph without multiple edges.
Further we assume that there are only finitely many directed paths between 
any two vertices.
If $u$ and $v$ are any pair of vertices in $D$,
let ${\cal P}(u,v)$ denote the set of all directed paths from $u$ to $v$ in $D$.
Fix a positive integer $n$.
An \defterm{$n$-vertex} is an $n$-tuple $\pmb{v}=(v_1,\dots,v_n)$ of $n$ vertices of $D$.
Given any pair of $n$-vertices $\pmb{u}=(u_1,\dots,u_n)$ and $\pmb{v}=(v_1,\dots,v_n)$,
an \defterm{$n$-path} from $\pmb{u}$ to $\pmb{v}$ is an $n$-tuple $\pmb{P}=(P_1,\dots,P_n)$ of $n$ paths such that $P_i\in{\cal P}(u_i,v_i)$.
Let ${\cal P}(\pmb{u},\pmb{v})$ denote the set of all $n$-paths from $\pmb{u}$ to $\pmb{v}$.
Two directed paths $P$ and $Q$ will be said to be non-intersecting if they share no common vertex.
An $n$-path $\pmb{P}$ is said to be non-intersecting if $P_i$ and $P_j$ are non-intersecting for any $i\not=j$.
Let ${\cal P}^0(\pmb{u},\pmb{v})$ denote the subset of ${\cal P}(\pmb{u},\pmb{v})$ which consists of all non-intersecting $n$-paths.

We assign a commutative indeterminate $x_e$ to each edge $e$ of $D$ and call it the \defterm{weight} of the edge.
Set the weight of a path $P$ to be the product of the weights of its edges and denote it by $\wt(P)$.
If $u$ and $v$ are any pair of vertices in $D$,
define
\begin{equation*}
h(u,v)=\sum_{P\in{\cal P}(u,v)} \wt(P).
\end{equation*}
The weight of an $n$-path is defined to be the product of the weights of its components.
The sum of the weights of $n$-paths in
${\cal P}(\pmb{u},\pmb{v})$ (resp. ${\cal P}^0(\pmb{u},\pmb{v})$) is denoted by $F(\pmb{u},\pmb{v})$ (resp. $F^0(\pmb{u},\pmb{v})$).
%
%
\begin{defi}
\label{dcompatible}
If $\pmb{u}=(u_1,\dots,u_n)$ and $\pmb{v}=(v_1,\dots,v_n)$ are $n$-vertices of $D$,
then $\pmb{u}$ is said to be \defterm{$D$-compatible} with $\pmb{v}$ if
every path $P\in{\cal P}(u_i,v_l)$ intersects with every path $Q\in{\cal P}(u_j,v_k)$ whenever $i<j$ and $k<l$.
\end{defi}
The following famous lemma is from \cite{GV}.
We recall its proof here again to make not only 
this paper self-contained
but also the subsequent discussion smooth.
%
%
\begin{lem}
\label{Gessel-Viennot}
(Lindstr\"om-Gessel-Viennot)
Let $\pmb{u}=(u_1,\dots,u_n)$ and $\pmb{v}=(v_1,\dots,v_n)$ be two $n$-vertices in an acyclic digraph $D$.
Then
\begin{equation}
\sum_{\pi\in S_n}\sgn\pi\ F^0(\pmb{u}^{\pi},\pmb{v})=\det[h(u_i,v_j)]_{1\le i,j\le n}.
\end{equation}
Here, for any permutation $\pi\in S_n$, let $\pmb{u}^{\pi}$ denote $(u_{\pi(1)},\dots,u_{\pi(n)})$.
In particular, if $\pmb{u}$ is $D$-compatible with $\pmb{v}$,
then
\begin{equation}
F^0(\pmb{u},\pmb{v})=\det[h(u_i,v_j)]_{1\le i,j\le n}.
\end{equation}
\end{lem}
%
%
{\it Proof.}
From the definition of determinants we have
\begin{equation}
\label{det_h}
\det[h(u_i,v_j)]_{1\le i,j\le n}=\sum_{\pi\in{\frak S}_n}
\sgn(\pi) h(u_1,v_{\pi(1)}) h(u_2,v_{\pi(2)}) \dots h(u_r,v_{\pi(n)}).
\end{equation}
For $\pi\in\frak{S}_n$, let $P(\pmb{u},\pmb{v}^\pi)$ denote the set of all the $n$-paths $\pmb{P}=\{P_1,\dots,P_n\}$ such that each path $P_i$ connects $u_i$ with $v_{\pi(i)}$ for $i=1,\dots,n$.
Let $P^0(\pmb{u},\pmb{v}^\pi)$ denote the subset of $P(\pmb{u},\pmb{v}^\pi)$ which consists of all non-intersecting paths $\pmb{P}\in P(\pmb{u},\pmb{v}^\pi)$.
Let us define sets $\Pi$ and $\Pi^0$ of configurations by
\begin{align*}
&\Pi=\left\{(\pi,\pmb{P}):\text{$\pi\in{\frak S}_n$ and
$\pmb{P}\in{\cal P}(\pmb{u},\pmb{v}^\pi))$}\right\},\\
&\Pi^0=\left\{(\pi,\pmb{P}):\text{$\pi\in{\frak S}_n$ and
$\pmb{P}\in{\cal P}^0(\pmb{u},\pmb{v}^\pi)$}\right\}.
\end{align*}
Then the right-hand side of \thetag{\ref{det_h}} is the generating function of configurations $(\pi,\pmb{P})\in\Pi$
with the weight $\wt(\pi,\pmb{P})=\sgn(\pi)\wt(\pmb{P})$.
Now we describe an involution on the set $\Pi\setminus\Pi^0$
which reverse the sign of the associated weight.
First fix an arbitrary total order on $V$.
Let $C=(\pi,\pmb{P})\in\Pi\setminus\Pi^0$.
Among all vertices that occurs as intersecting points,
let $v$ denote the least vertex with respect to the fixed order.
Among paths that pass through $v$,
assume that $P_i$ and $P_j$ are the two whose indices $i$ and $j$ are smallest.
Let $P_i(\rightarrow v)$ (resp. $P_i(v\rightarrow)$) denote the sub-path of $P_i$ from $u_i$ to $v$ (resp. from $v$ to $v_{\pi(i)}$).
Set $C'=(\pi',\pmb{P}')$ to be the configuration in which $P_k'=P_k$ for $k\not=i,j$,
\begin{equation*}
P_i'=P_i(\rightarrow v)P_j(v\rightarrow),\qquad
P_j'=P_j(\rightarrow v)P_i(v\rightarrow),
\end{equation*}
and $\pi'=\pi\circ(i,j)$.
It is easy to see that $C'\in\Pi$ and $\wt(C')=-\wt(C)$.
Thus $C\mapsto C'$ defines a sign reversing involution and, by this involution,
one may cancel all of the terms $\{\wt(C):C\in\Pi\setminus\Pi^0\}$
and only the terms $\{\wt(C):C\in\Pi^0\}$ remains.
Since $F^0(\pmb{u}^{\pi},\pmb{v})=F^0(\pmb{u},\pmb{v}^{\pi^{-1}})$,
we obtain the resulting identity.
In particular, if $\pmb{u}$ is $D$-compatible with $\pmb{v}$,
the configurations $C\in\Pi^0$ occur only when $\pi=\id$,
and are counted with the weight $\wt(P)$.
This proves the lemma.
$\Box$
\medbreak

\begin{figure}[b]
\begin{center}
\begin{picture}(200,60)
\put(  0,20){\circle*{8}}
\put( 30,20){\circle*{8}}
\put( 60,20){\circle*{6}}
\put( 90,20){\circle*{6}}
\put(120,20){\circle*{6}}
\put(150,20){\circle*{6}}
\put(180,20){\circle*{6}}
\put(210,20){\circle*{6}}
\put(  0,0){\makebox(0,0)[b]{$u_1$}}
\put( 30,0){\makebox(0,0)[b]{$u_2$}}
\put( 60,0){\makebox(0,0)[b]{$V_6$}}
\put( 90,0){\makebox(0,0)[b]{$V_5$}}
\put(120,0){\makebox(0,0)[b]{$V_4$}}
\put(150,0){\makebox(0,0)[b]{$V_3$}}
\put(180,0){\makebox(0,0)[b]{$V_2$}}
\put(210,0){\makebox(0,0)[b]{$V_1$}}
\qbezier(  0,20)( 60, 80)(120,20)
\qbezier( 30,20)(120,100)(210,20)
\qbezier( 60,20)(105, 60)(150,20)
\qbezier( 90,20)(135, 60)(180,20)
\end{picture}
\caption{Proof of Theorem~\ref{lmsf}}\label{figure:matching2}
\end{center}
\end{figure}
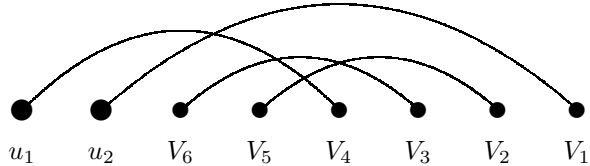

%
%
Let $n$ be an even integer and let $\pmb{v}=(v_1,\dots,v_n)$ be an $n$-vertex.
We write
\[
{\cal F}(\pmb{v})=\{(v_{\sigma_{1}},v_{\sigma_{2}},\dots,v_{\sigma_{n-1}},v_{\sigma_{n}})\,:\,
\sigma=(\sigma_{1},\sigma_{2},\hdots,\sigma_{n-1},\sigma_{n})\in{\cal F}_n\}
\]
and call it the set of perfect matchings of $\pmb{v}$.
Let $S=\{v_1<\cdots<v_{N}\}$ be a finite totally ordered subset of $V$.
We assume a commutative indeterminate $a_{v_iv_j}$ is assigned to each pair $(v_{i},v_{j})$ ($i<j$) of vertices in $S$.
We write the assembly of the indeterminates as $A=(a_{v_iv_j})_{i<j}$.
This upper triangular array uniquely defines a skew-symmetric matrix of size $N$,
and we use the same symbol $A$ to express this skew-symmetric matrix.
Suppose $m$ is even and $\pmb{u}=(u_1,\dots,u_m)$ is an $m$-vertex.
If $I=\{v_{i_1},\dots,v_{i_m}\}_{<}\in\binom{S}{m}$ is an $m$-element subset of $S$,
then we write $A^{I}_{I}$ for $(a_{v_{i_k}v_{i_k}})_{1\leq k<l\leq m}$.
If $\pmb{u}=(u_1,\dots,u_m)$ is an $m$-vertex and $\pmb{v}=(v_1,\dots,v_n)$ (resp. $S=\{v_1,\dots,v_n\}_{<}$) is an $n$-vertex (resp. an $n$-element totally ordered subset of $V$),
then let $H(\pmb{u},\pmb{v})$ (resp. $H(\pmb{u},S)$) denote the $m$ by $n$ matrix $\left(h(u_i,v_j)\right)_{1\leq i\leq v,\,1\leq j\leq n}$.
We consider the generating function of the set of non-intersecting $m$-paths from $\pmb{u}$ to $S$ weighted by subpfaffians of $A$:
\begin{equation*}
Q(\pmb{u},S;A)=\sum_{I\in\binom{S}{m}}
\Pf(A^{I}_{I})F^0(\pmb{u},I)
\end{equation*}
The following theorem express this generating function by a Pfaffian,
and is interpreted as a lattice path version of Theorem~\ref{msf}.
%
%
\begin{thm}
\label{lmsf}
Let $m$ and $N$ be even integers such that $0\leq m\leq N$.
Let $\pmb{u}=(u_1,\dots,u_m)$ be an $m$-vertex and $S=\{V_1<\cdots<V_N\}$ be a totally ordered set of vertices in an acyclic digraph $D$.
Let $A=(a_{V_iV_j})_{1\leq i<j\leq N}$ be an skew-symmetric matrix with rows and columns indexed by $S$,
and let $\widehat A$ denote its copfaffian matrix.
Then we have 
\begin{align}
\label{lmsf_eq_gen}
\sum_{I\in\binom{S}{m}}
\Pf(A^{I}_{I})\sum_{\pi\in S_m}\sgn\pi\,F^0(\pmb{u}^{\pi},I)
&=\Pf(A)\Pf\begin{pmatrix}
O_{m}& H(\pmb{u},S) J_N\\
- J_N\, {}^t H(\pmb{u},S)&\frac{1}{\Pf(A)}J_N {}^t\!\widehat A J_N\\
\end{pmatrix}.
\end{align}
In particular,
if $\pmb{u}$ is $D$-compatible with $S$,
then 
\begin{align}
\label{lmsf_eq}
\sum_{{I\subseteq S}\atop{\sharp I=m}}
\Pf(A^{I}_{I})F^0(\pmb{u}^{\pi},I)
&=\Pf(A)\Pf\begin{pmatrix}
O_{m}& H(\pmb{u},S) J_N\\
- J_N\, {}^t H(\pmb{u},S)&\frac{1}{\Pf(A)}J_N {}^t\!\widehat A J_N\\
\end{pmatrix}.
\end{align}
\end{thm}
%
%
{\it Proof.}
Let ${\widehat a}_{V_iV_j}$ denote the $(i,j)$-copfaffian of $A$
(i.e. $\widehat A=({\widehat a}_{V_iV_j})$)
and put $\alpha_{V_iV_j}=\frac1{\Pf(A)}\widehat a_{V_iV_j}$ for $1\leq i,j\leq N$.
Since multiplying $J_N$ from right reverse the order of columns of $H(\pmb{u},S)$,
we have
\begin{align*}
\Pf\begin{pmatrix}
O_{n}& H(\pmb{u},S) J_N\\
- J_N\, {}^t H(\pmb{u},S)&\frac{1}{\Pf(A)}J_N {}^t\!\widehat A J_N\\
\end{pmatrix}
=\sum_{\tau}\sgn\tau\prod_{(u_i,V_j)\in\tau}h(u_i,V_j)\prod_{(V_k,V_l)\in\tau}\alpha_{V_kV_l}
\end{align*}
summed over all perfect matchings $\tau$ on $(u_1,\dots,u_m,V_N,\dots,V_1)$
in which there are no edges connecting any two vertices of $u$.	
For an example of such a perfect matching, see Figure~\ref{figure:matching2} bellow.
We may interpret this Pfaffian
as the generating function for all $(m+1)$-tuples $C=(\tau,P_1,\dots,P_m)$
such that $P_i\in{\cal P}(u_i,V_j)$ if there is an edge $(u_i,V_j)\in\tau$.
This implies that each vertex, say $u_i$, in $\pmb{u}$ is always connected to a vertex,
say $V_j$, in $S$,
and remaining $(N-n)$ vertices of $S$ are connected each other by edges
which we write $(V_k,V_l)\in\tau$.
The weight assigned to $C=(\tau,P_1,\dots,P_m)$	shall be 
$\sgn\tau\,\left(\prod_{(V_k,V_l)\in\tau}\alpha_{V_kV_l}\right) w(P_1)\cdots w(P_m)$.
Let $\Sigma$ denote the set of configurations $C=(\tau,P_1,\dots,P_m)$ satisfying the above condition,
and let $\Sigma^0$ denote the subset consisting of all configurations $C=(\tau,P_1,\dots,P_m)$ such that
$(P_1,\dots,P_m)$ is non-intersecting.
We will show that there is a sign-reversing involution on $\Sigma\setminus\Sigma^0$,
i.e. the set of the configurations $C=(\tau,P_1,\dots,P_m)$ with at least one pair of intersecting paths.
Our proof here is essentially the same as that in Lemma~\ref{Gessel-Viennot}.
To describe the involution,
first choose a fixed total order of the vertices,
and consider an arbitrary configuration $C=(\tau,P_1,\dots,P_m)\in\Sigma\setminus\Sigma^0$.
Among all vertices that occurs as intersecting points,
let $v$ denote the vertex which precedes all other points of intersections with respect to the fixed order.
Among paths that pass through $v$,
assume that $P_i$ and $P_j$ are the two whose indices $i$ and $j$ are smallest.
We define $C'=(\tau',P_1',\dots,P_m')$ to be the configuration where
$P_{k}'=P_{k}$ for $k\not=i,j$,
\begin{equation*}
P_i'=P_i(\rightarrow v)P_j(v\rightarrow),\qquad
P_j'=P_j(\rightarrow v)P_i(v\rightarrow),
\end{equation*}
and, if $(u_i,V_k)$ and $(u_j,V_l)$ are the edges of $\tau$,
then $(u_i,V_l)$ and $(u_j,V_k)$ are in $\tau'$ and all the other edges of $\tau'$ are the same as $\tau$.
Note that the multi-sets of edges appearing in $(P_1,\dots,P_m)$ and $(P_1',\dots,P_m')$ are identical,
which means $\wt(C')=-\wt(C)$.
Since this involution changes the sign of the associated weight,
one may cancel all of the terms appearing in $\Sigma\setminus\Sigma^0$,
aside from those with non-intersecting paths.
For $C=(\tau,P_1,\dots,P_m)\in\Sigma^0$,
let $I$ denote the set of vertices of $S$ connected to a vertex in $\pmb{u}$,
and let $\overline I$ denote the complementary set of $I$ in $S$.
Put $I=\{V_{i_1},\dots,V_{i_m}\}_{<}$,
then we can find a unique permutation $\pi\in S_m$ such that each $u_{\pi(k)}$ is connected to $V_{i_{k}}$ in $\tau$ for $k=1,\dots,m$.
The remaining edges in $\tau$ which does not contribute to this permutation
perform a perfect matching on $\overline I$ which we denote by $\sigma$.
The $\sgn\tau$ is equal to $(-1)^{s(I,\overline I)}\sgn\pi\sgn\sigma$,
where $s(I,\overline I)$ denote the shuffle number to merge $I$ with $\overline I$ into $S$.
Thus,
if we put $m=2m'$ and $N=2N'$ for nonnegative integers $m'$ and $N'$,
the sum of weights is equal to
\begin{equation*}
\sum_{I}(-1)^{s(I,\overline I)}\sum_{\pi\in S_m}\sgn\pi\,F^{0}(\pmb{u}^{\pi},I)\frac{1}{(\Pf A)^{N'-m'}}\Pf \left(\widehat A_{\overline I}^{\overline I}\right),
\end{equation*}
where $I$ runs over all subsets of $S$ of cardinality $m$.
Theorem~\ref{Lewis-Carroll_pfaffian}
implies
\begin{equation*}
\Pf \left(\widehat A_{\overline I}^{\overline I}\right)=(-1)^{|\overline I|+N'-m'}\Pf(A)^{N'-m'-1}\Pf(A_{I}^{I}).
\end{equation*}
Since $I\cup\overline I=S$,
we have $|I|+|\overline I|=\binom{N+1}{2}\equiv N'$ ($\module2$).
Meanwhile, it is easy to see $(-1)^{s(I,\overline I)}=(-1)^{|I|-m'}$.
This immediately implies \thetag{\ref{lmsf_eq_gen}}.
This completes the proof.
$\Box$
\medbreak

In fact Theorem~\ref{lmsf} is equivalent to Theorem~\ref{msf}.
In \cite{IW1} we gave an algebraic proof of Theorem~\ref{msf} using the exterior algebra.
One clearly sees that Theorem~\ref{lmsf} is an easy consequence of Theorem~\ref{msf} and Lemma~\ref{Gessel-Viennot}.
Here we give a proof that derives Theorem~\ref{msf} from Theorem~\ref{lmsf}.
Similarly one can derive Theorem~\ref{msf2} from Theorem~\ref{lmsf2}, and also Theorem~\ref{msf3} from Theorem~\ref{lmsf3},
but we will not give the details here and leave it to the reader.

\medbreak
\noindent
%
%
{\it Proof of Theorem~\ref{msf}}.
First we define a digraph $D$ with vertex set $\Bbb{Z}^2$ and edges directed from $u$ to $v$ whenever $v-u=(1,0)$ or $(0,1)$.
For $u=(i,j)$, we assign the weight $x_j$ (resp. $1$) to the edge with $v-u=(1,0)$ (resp. $(0,1)$).
If $u=(i,1)$ and $v=(j,r)$, then $\lim_{r\rightarrow\infty}h(u,v)=h_{j-i}(x)$ is well-known as the complete symmetric function,
which is defined by $\sum_{k\ge0}h_k(x)t^k=\prod_{i\ge1}\frac1{1-x_it}$
(See \cite{Ma}).
Thus, if we fix constants $(\lambda_1,\cdots,\lambda_{m})$ and $(\mu_1,\cdots,\mu_{m})$
which satisfy $\lambda_1<\cdots<\lambda_m$ and $\mu_1<\cdots<\mu_m$,
and take the vertices $u_i=(\lambda_i,1)$ and $v_i=(\mu_i,r)$ for $i=1,\cdots,m$,
then, $\pmb{u}$ and $\pmb{v}$ are $D$-compatible,
and from Lemma~\ref{Gessel-Viennot}, we deduce
\begin{equation*}
\lim_{r\rightarrow\infty}F^0(\pmb{u},\pmb{v})=\det\left(h_{\mu_j-\lambda_i}(x)\right)_{1\le i,j\le m}.
\end{equation*}
Let $N$ be a positive integer such that $N\ge m\geq0$ and $A=(a_{ij})$ be an $N$ by $N$ skew-symmetric matrix.
We let $\pmb{u}=(u_1,\cdots,u_m)$ with $u_i=(Ni,1)$ for $i=1\cdots m$ and $S=\{v_1,\cdots,v_N\}$ with $v_j=(j+Nm,r)$ for $j=1\cdots N$.
Then, from Theorem~\ref{lmsf} and by taking the limit $r\rightarrow\infty$, we obtain
\begin{equation*}
\sum_{I=\{i_1<\cdots<i_m\}\subseteq[N]}
\Pf(A_I^I)\det\left(h_{i_j+N(m-i)}(x)\right)_{1\le i,k\le m}
=\Pf\left(Q\right)
\end{equation*}
where $Q=(Q_{ij})_{1\le i,j\le m}$ is given by
\begin{equation*}
Q_{ij}=\sum_{1\le k<l\le N} a_{kl}
\begin{vmatrix}
h_{k+N(n-i)}(x)&h_{l+N(n-i)}(x)\\
h_{k+N(n-j)}(x)&h_{l+N(n-j)}(x)\\
\end{vmatrix}
\end{equation*}
We use the fact that the $h_1$, $\cdots$, $h_{Nm}$ are algebraically independent over $\Bbb{Q}$
(See \cite{Ma}).
Thus we can replace each $h_{j+N(m-i)}$ with any commutative indeterminate $t_{ij}$,
and we obtain Theorem~\ref{msf}.
$\Box$
\bigbreak

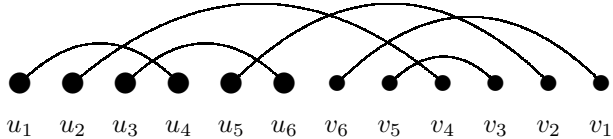
\begin{figure}[b]
\begin{center}
\begin{picture}(200,60)
\put(  0,20){\circle*{8}}
\put( 20,20){\circle*{8}}
\put( 40,20){\circle*{8}}
\put( 60,20){\circle*{8}}
\put( 80,20){\circle*{8}}
\put(100,20){\circle*{8}}
\put(120,20){\circle*{6}}
\put(140,20){\circle*{6}}
\put(160,20){\circle*{6}}
\put(180,20){\circle*{6}}
\put(200,20){\circle*{6}}
\put(220,20){\circle*{6}}
\put(  0,0){\makebox(0,0)[b]{$u_1$}}
\put( 20,0){\makebox(0,0)[b]{$u_2$}}
\put( 40,0){\makebox(0,0)[b]{$u_3$}}
\put( 60,0){\makebox(0,0)[b]{$u_4$}}
\put( 80,0){\makebox(0,0)[b]{$u_5$}}
\put(100,0){\makebox(0,0)[b]{$u_6$}}
\put(120,0){\makebox(0,0)[b]{$v_6$}}
\put(140,0){\makebox(0,0)[b]{$v_5$}}
\put(160,0){\makebox(0,0)[b]{$v_4$}}
\put(180,0){\makebox(0,0)[b]{$v_3$}}
\put(200,0){\makebox(0,0)[b]{$v_2$}}
\put(220,0){\makebox(0,0)[b]{$v_1$}}
\qbezier(  0,20)( 30, 50)( 60,20)
\qbezier( 20,20)( 90, 80)(160,20)
\qbezier( 40,20)( 70, 50)(100,20)
\qbezier(120,20)(170, 70)(220,20)
\qbezier( 80,20)(140, 80)(200,20)
\qbezier(140,20)(160, 40)(180,20)
\end{picture}
\caption{Proof of Theorem~\ref{lmsf3}}\label{figure:matching4}
\end{center}
\end{figure}

We next consider the lattice path version of Theorem~\ref{msf2}.
Let $m$, $n$ and $N$ be positive integers such that $m-n$ is even and $0\leq m-n\leq N$.
Suppose that $S^0=\{v_1<\dots<v_n\}$ is a fixed ordered list of vertices,
and let $S=\{V_{1}<\cdots<V_{N}\}$ be a totally ordered set disjoint with $S^0$.
For a subset $I$ of $S$
let $S^0\uplus I$ denote the union of $S^0$ and $I$,
ordered so that each $v_i$ precedes each $w\in I$.
Let $A=(a_{V_iV_j})_{1\leq i,j\leq N}$ be a skew-symmetric matrix with rows and columns indexed by the totally ordered set $S$ as before.
We will obtain a formula of the generating function weighted by subpfaffians of $A$ as follows:
\begin{equation*}
Q(\pmb{u};S^0,S;A)=\sum_{I\in\binom{S}{m-n}}
\Pf(A_I^I)F^0(\pmb{u},S^0\uplus I)
\end{equation*}
\medbreak
%
%
\begin{thm}
\label{lmsf2}
Let $m$, $n$ and $N$ be positive integers such that $m-n$ and $N$ are even integers and $0\leq m-n\leq N$.
Let $\pmb{u}=(u_1,\cdots,u_m)$ be an $m$-vertex and $S^0=\{v_1<\dots<v_n\}$ be an $n$-vertex in an acyclic digraph $D$.
Let $S=\{V_{1}<\cdots<V_{N}\}$ be a finite totally ordered set of vertices which is disjoint with $S_0$.
Let $A=(a_{V_iV_j})$ be an $N$ by $N$ skew-symmetric matrix with rows and columns indexed by $S$.
Then
\begin{align}
&\sum_{I\in\binom{S}{m-n}}
\Pf(A^{I}_{I})\sum_{\pi\in S_m}\sgn\!\pi\,\pmb{F}^0(u^{\pi},S^0\uplus I)
\nonumber\\
&=\Pf(A)
\Pf\begin{pmatrix} 
   O_{m}   &   H(\pmb{u};S) J_N & H(\pmb{u};S^0) J_n \\
- J_N\, {}^tH(\pmb{u};S) & \frac1{\Pf(A)}J_N {}^t\!\widehat A J_N& O_{N,n}\\
- J_n\, {}^tH(\pmb{u};S^0) & O_{n,N} & O_{n}
\end{pmatrix}
\label{lmsf2_eq_gen}
\end{align}
where
\begin{align*}
&H(\pmb{u};S^0)=\left(h(u_i,v_j)\right)_{1\le i\le m, 1\le j\le n},\\
&H(\pmb{u};S)=\left(h(u_i,V_j)\right)_{1\le i\le m, 1\le j\le N}.
\end{align*}
In particular, if $\pmb{u}$ is D-compatible with $I^0\uplus S$,
then we have
\begin{align}
&\sum_{I\in\binom{S}{m-n}}
\Pf(A^{I}_{I})\pmb{F}^0(u,S^0\uplus I)
\nonumber\\
&=\Pf(A)
\Pf\begin{pmatrix} 
   O_{m}   &   H(\pmb{u};S) J_N & H(\pmb{u};S^0) J_n \\
- J_N\, {}^tH(\pmb{u};S) & \frac1{\Pf(A)}J_N {}^t\!\widehat A J_N& O_{N,n}\\
- J_n\, {}^tH(\pmb{u};S^0) & O_{n,N} & O_{n}
\end{pmatrix}
\label{lmsf2_eq}
\end{align}
\end{thm}
%
%
{\it Proof}.
Let $\widehat a_{V_kV_l}$ denote the $(k,l)$-copfaffain of $A$,
and put $\alpha_{V_kV_l}=\frac1{\Pf(A)}\widehat a_{V_kV_l}$ as before.
We have
\begin{align}
&\Pf\begin{pmatrix} 
   O_{m}   &   H(\pmb{u};S) J_N & H(\pmb{u};S^0) J_n \\
- J_N {}^tH(\pmb{u};S) & \frac1{\Pf(A)}J_N {}^t\!\widehat A J_N& O_{N,n}\\
- J_n {}^tH(\pmb{u};S^0) & O_{n,N} & O_{n}
\end{pmatrix}\nonumber\\
&=\sum_{\tau} \sgn\tau \prod_{(u_i,V_j)\in\tau} h(u_i,V_j) \prod_{(V_i,V_j)\in\tau} \alpha_{V_iV_j} \prod_{(u_i,V_j)\in\tau} h(u_i,V_j)
\label{pf_interpret}
\end{align}
summed over all perfect matchings $\tau$ of $(u_1,\dots,u_m,V_N,\dots,V_1,v_n,\dots,v_1)$
in which there are no edges connecting any two vertices of $\pmb{u}$,
and each vertex in $S^0$ must be connected to a vertex in $\pmb{u}$.
An example of such a perfect matching is given below.
\begin{figure}[b]
\begin{center}
\begin{picture}(200,60)
\put(  0,20){\circle*{8}}
\put( 20,20){\circle*{8}}
\put( 40,20){\circle*{8}}
\put( 60,20){\circle*{8}}
\put( 80,20){\circle*{6}}
\put(100,20){\circle*{6}}
\put(120,20){\circle*{6}}
\put(140,20){\circle*{6}}
\put(160,20){\circle*{6}}
\put(180,20){\circle*{6}}
\put(200,20){\circle*{8}}
\put(220,20){\circle*{8}}
\put(  0,0){\makebox(0,0)[b]{$u_1$}}
\put( 20,0){\makebox(0,0)[b]{$u_2$}}
\put( 40,0){\makebox(0,0)[b]{$u_3$}}
\put( 60,0){\makebox(0,0)[b]{$u_4$}}
\put( 80,0){\makebox(0,0)[b]{$V_6$}}
\put(100,0){\makebox(0,0)[b]{$V_5$}}
\put(120,0){\makebox(0,0)[b]{$V_4$}}
\put(140,0){\makebox(0,0)[b]{$V_3$}}
\put(160,0){\makebox(0,0)[b]{$V_2$}}
\put(180,0){\makebox(0,0)[b]{$V_1$}}
\put(200,0){\makebox(0,0)[b]{$v_2$}}
\put(220,0){\makebox(0,0)[b]{$v_1$}}
\qbezier(  0,20)(100,100)(200,20)
\qbezier( 20,20)( 70, 60)(120,20)
\qbezier( 40,20)(130,100)(220,20)
\qbezier( 60,20)(120, 60)(160,20)
\qbezier( 80,20)(130, 60)(180,20)
\qbezier(100,20)(120, 40)(140,20)
\end{picture}
\caption{Proof of Theorem~\ref{lmsf2}}\label{figure:matching3}
\end{center}
\end{figure}
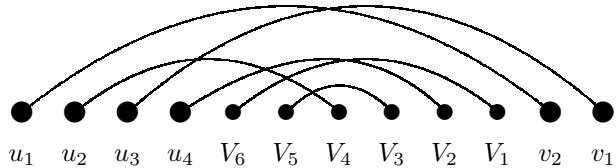
This may be interpreted as the generating function for all $(m+1)$-tuples
$C=(\tau,P_1,\dots,P_m)$ such that $P_i\in{\cal P}(u_i,v_j)$ if there is an edge $(u_i,v_j)\in\tau$,
and $P_i\in{\cal P}(u_i,V_j)$ if there is an edge $(u_i,V_j)\in\tau$.
The weight assigned to $C=(\tau,P_1,\dots,P_m)$ shall be 
$\sgn\tau\,\prod_{(V_k,V_l)\in\tau}\alpha_{V_kV_l}w(P_1)\dots w(P_m)$.
We claim that the sign-reversing involution used in the previous proofs can be applied to this situation as well.
In fact, quite the same arguments show that one may cancel all of the terms appearing in \thetag{\ref{pf_interpret}},
aside from those with non-intersecting paths.
In $\tau$ associated with these configuration $C=(\tau,P_1,\dots,P_m)$,
each $v_k$ ($k=1,\dots,n$) is always connected to a vertex in $\pmb{u}$.
This means that exactly $n$ vertices of $\pmb{u}$ are connected to vertices in $S^0$,
and the remaining $(m-n)$ vertices are connected to certain vertices in $S$.
Let $I$ denote the set of vertices in $S$ connected to vertices in $\pmb{u}$,
and let $\overline I$ denote its complementary set in $S$.
Note that $\sharp I=(m-n)$ and $\sharp\overline I=(N-m+n)$.
Let $S^0\uplus I$ denote the juxtaposition of vertices from $I^0$ and $I$ arranged in this order,
and, if we put $I^0\uplus I=(u_1^{*},\dots,u_m^{*})$
then there is a unique permutation $\pi\in S_{m}$ such that each $u_{\pi(k)}$ is connected to $u_{k}^{*}$ for $k=1,\dots,m$.
The remaining edges of $\tau$ whose both endpoints are included $\overline I$ define a unique perfect matching on $\overline I$.
In this situation $\sgn\tau$ is equal to $s(I,\overline I)\sgn\pi\,\sgn\sigma$.
Thus,
if we put $m-n=2m'$ and $N=2N'$ for nonnegative integers $m'$ and $N'$,
then the sum of weights is equal to
\begin{equation*}
\sum_{I}(-1)^{s(I,\overline I)}\sgn\pi\,\pmb{F}^0(u^{\pi},I^0\uplus I)\frac1{(\Pf A)^{N'-m'}}\Pf({\widehat A}^{\overline I}_{\overline I}),
\end{equation*}
where $I$ runs over all subsets of $S$ of cardinality $(m-n)$.
From Theorem~\ref{Lewis Carroll_pfaffian}
we have
\begin{equation*}
\Pf \left(\widehat A_{\overline I}^{\overline I}\right)=(-1)^{|\overline I|+N'-m'}(\Pf A)^{N'-m'-1}\Pf(A_{I}^{I}).
\end{equation*}
Since $I\cup\overline I=S$,
we have $|I|+|\overline I|=\binom{N+1}{2}\equiv N'$ ($\module2$).
Meanwhile, it is easy to see $(-1)^{s(I,\overline I)}=(-1)^{|I|-m'}$.
This immediately implies \thetag{\ref{lmsf2_eq_gen}}.
Lastly, if $\pmb{u}$ is D-compatible with $I^0\uplus S$,
then there is no non-intersecting path unless $\pi=id$,
which immediately implies \thetag{\ref{lmsf2_eq}}.
This completes the proof.
$\Box$
\medbreak
%
%
Now we show the following theorem for proving Theorem~\ref{msf3}.
For the purpose we consider a more general problem concerning with
non-intersecting paths in which both starting points and end points vary. 
%
%
\begin{thm}
\label{lmsf3}
Let $M$ and $N$ be even integers.
Let $R=\{u_1<\dots<u_M\}$ and $S=\{v_1<\dots<v_N\}$ be totally ordered subsets of vertices in an acyclic digraph $D$.
Let $A=(a_{u_iu_j})$ (resp. $B=(b_{v_iv_j}$) be a non-singular skew-symmetric matrix with rows and columns indexed by the vertices of $R$ (resp. $S$).
Then
\begin{align}
&\sum_{{0\leq r\leq \min(M,N)}\atop{\text{$r$ even}}} z^{r}
\sum_{I\in\binom{R}{r}}
\sum_{J\in\binom{S}{r}} \Pf(A_I^I)\Pf(B^J_J)
\sum_{\pi\in S_n}\sgn\pi\,F^0(I^{\pi},J)
\nonumber\\
&=\Pf(A)\Pf(B)\Pf\begin{pmatrix}
\frac1{\Pf(A)}\widehat A& zH(R,S) J_N\\
- zJ_N\, {}^t H(R,S)&\frac{1}{\Pf(B)}J_N {}^t\!\widehat B J_N\\
\end{pmatrix}.
\label{row_column_even}
\end{align}
In particular,
if $R$ is compatible with $S$,
then
\begin{align}
&\sum_{{0\leq r\leq\min(M,N)}\atop{\text{$r$ even}}} z^{r}
\sum_{I\in\binom{R}{r}}
\sum_{J\in\binom{S}{r}}
\Pf(A_I^I)\Pf(B^J_J)
F^0(I,J)
\nonumber\\
&=\Pf(A)\Pf(B)\Pf\begin{pmatrix}
\frac1{\Pf(A)}\widehat A& zH(R,S) J_N\\
- zJ_N\, {}^t H(R,S)&\frac{1}{\Pf(B)}J_N {}^t\!\widehat B J_N\\
\end{pmatrix}.
\label{row_column_even2}
\end{align}
\end{thm}
%
%
{\it Proof}.
Let $\widehat a_{u_iu_j}$ (resp. $\widehat b_{v_iv_j}$) denote the $(i,j)$-copfaffain of $A$ (resp. $B$).
Put $\alpha_{u_iu_j}=\frac1{\Pf(A)}\widehat a_{u_iu_j}$ and $\beta_{v_iv_j}=\frac1{\Pf(B)}\widehat b_{v_iv_j}$.
Then, we have
\begin{align*}
&\Pf\begin{pmatrix}
\frac1{\Pf(A)}\widehat A& zH(R,S) J_N\\
- zJ_N\, {}^t H(R,S)&\frac{1}{\Pf(B)}J_N {}^t\!\widehat B J_N\\
\end{pmatrix}\\
&=\sum_{\tau}\sgn\tau\prod_{(u_i,u_j)\in\tau}\alpha_{u_iu_j} \prod_{(u_i,v_j)\in\tau}z h(u_i,v_j)\prod_{(v_i,v_j)\in\tau}\beta_{v_iv_j}
\end{align*}
summed over all perfect matchings $\tau$ on $(u_1,\dots,u_M,v_N,\dots,v_1)$.
An example of such a perfect matching is Figure~\ref{figure:matching4} bellow.
As before we may interpret this Pfaffian as the generating function for all $(r+1)$-tuples
$C=(\tau,P_1,\dots,P_r)$ which satisfies (i) $r$ is an even integer such that $0\leq r\leq \min(M,N)$,
(ii) there are exactly $r$ edges whose one endpoint is in $R$ and the other endpoint is $S$,
and (iii) $\tau$ is a perfect matching on $(u_1,\dots,u_M,v_N,\dots,v_1)$ such that $P_i\in{\cal{P}}(u_i,v_j)$ if and only if there is an edge $(u_i,v_j)\in\tau$.
The weight assigned to $C=(\tau,P_1,\dots,P_r)$ shall be 
\begin{equation*}
\sgn\tau\:z^r\prod_{(u_i,u_j)\in\tau}\alpha_{u_iu_j}\:w(P_1)\cdots w(P_r)\prod_{(v_i,v_j)\in\tau}\beta_{v_iv_j}.
\end{equation*}
The same argument as in the proof of Theorem~\ref{lmsf} shows us that
we can define a sign-reversing involution on the set of the configurations $C=(\tau,P_1,\dots,P_r)$
with at least one pair of intersecting paths,
and this involution cancels all of the terms involving intersecting configurations of paths.
Thus we need to sum over only non-intersecting configurations.
Given a perfect matching $\tau$ on $(u_1,\dots,u_M,v_N,\dots,v_1)$
such that there are exactly $r$ edges connecting a vertex in $R$ and a vertex in $S$.
Let $I$ (resp. $J$) denote the subset of $R$ (resp. $S$)
which is composed of such endpoints of $\tau$.
Thus $\sharp I=\sharp J=r$, and $r$ must be even.
Let $\overline I$ (resp. $\overline J$) denote the complementary set of $I$ (resp. $J$) in $R$ (resp. $S$).
Put $I=\{u_{i_1},\dots,u_{i_r}\}_{<}$ and $J=\{v_{j_1},\dots,v_{j_r}\}_{<}$,
then there is a unique permutation $\pi$ such that $u_{i_{\pi(\nu)}}$ is connected to $v_{j_\nu}$ in $\tau$ for $\nu=1,\dots,r$.
If we put $M=2M'$, $N=2N'$ and $r=2r'$,
then the sum of weights becomes
\begin{equation*}
\sum_{{r=0}\atop{r=2r'}}^{M}
\sum_{{I\subseteq R}\atop{\sharp I=r}}\sum_{{J\subseteq S}\atop{\sharp J=r}}
(-1)^{s(I,\overline I)+s(J,\overline J)}
\sum_{\pi\in S_n}\sgn\pi\,z^rF^0(I^{\pi},J)
\frac1{\Pf(A)^{M'-r'}}\Pf({\widehat A}^{\overline I}_{\overline I})
\frac1{\Pf(B)^{N'-r'}}\Pf({\widehat B}^{\overline J}_{\overline J}).
\end{equation*}
By \thetag{\ref{LCpfaffian}} we have 
$\Pf({\widehat A}^{\overline I}_{\overline I})
=(-1)^{|\overline I|-M'+r'}\Pf(A)^{M'-r'-1}\Pf(A^I_I)$
and
$\Pf({\widehat B}^{\overline J}_{\overline J})
=(-1)^{|J|-N'+r'}\Pf(B)^{N'-r'-1}\Pf(B^J_J)$.
Further it is easy to see that $s(I,\overline I)\equiv |I|-r'$ ($\MOD2$) and $s(J,\overline J)\equiv |J|-r'$ ($\MOD2$).
Since $I\cup\overline I=R$ and $J\cup\overline J=S$,
we have $|I|+|\overline I|=\binom{M+1}{2}\equiv M'$ ($\MOD2$) 
and $|J|+|\overline J|=\binom{N+1}{2}\equiv N'$ ($\MOD2$).
These identities immediately implies \thetag{\ref{row_column_even}}.
$\Box$

\medbreak
In the following theorem we assume $D$-compatibility of two regions to make our notation simple,
but a more general theorem is also possible to establish.
%
%
\begin{thm}
\label{lmsf4}
Let $m$, $n$, $M$ and $N$ be nonnegative integers such that $M\equiv N\equiv m-n\equiv0$ ($\MOD2$).
Let $R^0=(u_1,\dots,u_m)$ (resp. $S^0=(v_1,\dots,v_n)$) be an $m$-vertex (resp. an $n$-vertex) in an acyclic digraph $D$.
Let $R=\{U_1<\dots<U_M\}$ (resp. $S=\{V_1<\dots<V_N\}$) be totally ordered subsets of vertices in $D$
which is disjoint with $R^0$ (resp. $S^0$).
Assume that $R^0\uplus R$ is $D$-compatible with $S^0\uplus S$.
Let $A=(a_{U_iU_j})$ (resp. $B=(b_{V_iV_j})$) be a skew-symmetric matrix with rows and columns indexed by the vertices of $R$ (resp. $S$).
Then
\begin{align}
&\sum_{{\max(m,n)\leq r\leq\min(m+M,n+N)}\atop{\text{$r-\max(m,n)$ is even.}}} z^{r}
\sum_{I\in\binom{R}{r-m}}
\sum_{J\in\binom{S}{r-n}}
\Pf(A_I^I)\Pf(B^J_J)
F^0(R^0\uplus I,S^0\uplus J)
\nonumber\\
&=\Pf(A)\Pf(B)
\Pf\begin{pmatrix}
O_m                           &  O_{m,M}  &  zH(R^0,S)J_N  & zH(R^0,S^0) J_n\\
O_{M,m}                       &\frac1{\Pf(A)}\widehat A   & zH(R;S) J_N & zH(R;S^0) J_n \\
-zJ_N {}^t\! H(R^0,S)   &-J_N{}^t\! zH(R;S) & \frac1{\Pf(B)}J_N {}^t\!\widehat B J_N& O_{N,n}\\
-zJ_n  {}^t\! H(R^0,S^0)&-J_n{}^t\! zH(R;S^0) & O_{n,N} & O_{n}
\end{pmatrix}
\end{align}
\end{thm}
%
%
{\it Proof}.
Without loss of generality we may assume that $m\geq n$.
Let $\widehat a_{U_iU_j}$ (resp. $\widehat b_{V_iV_j}$) denote the $(i,j)$-copfaffain of $A$ (resp. $B$).
Put $\alpha_{U_iU_j}=\frac1{\Pf(A)}\widehat a_{U_iU_j}$ and $\beta_{V_iV_j}=\frac1{\Pf(B)}\widehat b_{V_iV_j}$.
Further we put $R^0\uplus R=(u_1,\dots,u_m,U_1,\dots,U_M)=(u_1^{*},\dots,u_{m+M}^{*})$ and $S^0\uplus S=(v_1,\dots,v_n,V_1,\dots,V_N)=(v_1^{*},\dots,v_{n+N}^{*})$ for convenience.
Then we have
\begin{align}
&\Pf\begin{pmatrix}
O_m                           &  O_{m,M}  &  zH(R^0,S)J_N  & zH(R^0,S^0) J_n\\
O_{M,m}                       &\frac1{\Pf(A)}\widehat A   & zH(R;S) J_N & zH(R;S^0) J_n \\
-zJ_N {}^t\! H(R^0,S)   &-zJ_N{}^t\! H(R;S) & \frac1{\Pf(B)}J_N {}^t\!\widehat B J_N& O_{N,n}\\
-zJ_n  {}^t\! H(R^0,S^0)&-zJ_n{}^t\! H(R;S^0) & O_{n,N} & O_{n}
\end{pmatrix}
\nonumber\\
&=\sum_{\tau}\sgn\tau\,\prod_{(U_k,U_l)\in\tau}\alpha_{U_kU_l}\prod_{(V_k,V_l)\in\tau}\beta_{V_kV_l}
\prod_{(u_k^{*},v_l^{*})\in\tau} z h(u_k^{*},v_l^{*})
\label{gf4}
\end{align}
summed over all perfect matching on $(u_1^{*},\dots,u_{m+M}^{*},v_{n+N}^{*},\dots,v_1^{*})$.
We can interpret this Pfaffian as the generating function for all $(r+1)$-tuples $C=(\tau,P_1,\dots,P_{r})$
which satisfies 
(i) $r$ is an integer such that $m\leq r\leq \min(m+M,n+N)$ and $r\equiv m$ ($\MOD2$),
(ii) $\tau$ is a perfect matching on $(u_1^{*},\dots,u_{m+M}^{*},v_{n+N}^{*},\dots,v_1^{*})=(u_1,\dots,u_m,U_1,\dots,U_M,V_N,\dots,V_1,v_n,\dots,v_1)$ 
such that $P_i\in {\cal P}(u_k^{*},v_l^{*})$ if and only if there is an edge $(u_k^{*},v_l^{*})\in\tau$ which is connecting a vertex in $R^0\uplus R$ and a vertex in $S^0\uplus S$.
(iii) each vertex in $R^0$ must be connected to a vertex in $S^0\uplus S$,
(iv) each vertex in $S^0$ must be connected to a vertex in $R^0\uplus R$,
(v) and there are exactly $r$ edges connecting a vertex in $R^0\uplus R$ with a vertex in $S^0\uplus S$.
The weight assigned to $C=(\tau,P_1,\dots,P_{r})$ shall be $\sgn\tau\,z^r\prod_{(U_k,U_l)\in\tau}\alpha_{U_kU_l}\prod_{(V_k,V_l)\in\tau}\beta_{V_kV_l}\wt(P_1)\cdots\wt(P_{r})$.
It is easy to see that the sign-reversing involution used in the previous proofs is applicable exactly as before,
and we may cancel all the terms appearing in \thetag{\ref{gf4}},
aside from those with non-intersecting paths.
Thus we only need to consider configurations $C=(\tau,P_1,\dots,P_{r})$ with non-intersecting paths.
From the assumption that $R^0\uplus R$ is $D$-compatible with $S^0\uplus S$,
(a) each $v_i$, $i=1,\dots,n$, must be connected to $u_i$ in $\tau$ and $P_i\in{\cal P}(u_i,v_i)$,
(b) there are $r-m$ vertices $I=\{U_{i_1},\dots,U_{i_{r-m}}\}_{1\leq i_1<\dots<i_{r-m}\leq M}$ in $R$ 
and $r-n$ vertices $J=\{V_{j_1},\dots,V_{j_{r-n}}\}_{1\leq j_1<\dots<j_{r-n}\leq N}$
such that 
each $u_k$, $k=n+1,\dots,m$, is connected to $V_{j_{k-n}}$ in $\tau$ and $P_k\in{\cal P}(u_k,V_{j_{k-n}})$
and each $U_{i_k}$, $k=1,\dots,r-m$, is connected to $V_{j_{k+m-n}}$ in $\tau$ and $P_{k+m}\in{\cal P}(U_{i_k},V_{j_{k+m-n}})$,
(d) and the remaining $(m+M-r)$ vertices in $R$ are connected each other in $\tau$ and the remaining $(n+N-r)$ vertices in $S$ are connected each other in $\tau$.
We set $\overline I$ (resp. $\overline J$) to be the complementary set of $I$ (resp. $J$) in $R$ (resp. $S$).
If we put $M=2M'$, $N=2N'$, $m-n=2l'$ and $r-m=2r'$ for nonnegative integers $M'$, $N'$, $l'$ and $r'$,
then the sum of the weights becomes
\begin{align*}
\sum_{r'=0}^{\min(M',N'-l')}\sum_{I\in\binom{R}{2r'}}\sum_{J\in\binom{S}{2(r'+l')}}
(-1)^{s(I,\overline I)+s(J,\overline J)}\frac1{\Pf(A)^{M'-r'}\Pf(B)^{N'-r'-l'}}
\Pf({\widehat A}^{\overline I}_{\overline I})\Pf({\widehat B}^{\overline J}_{\overline J})F^0(R^0\uplus I,S^0\uplus J).
\end{align*}
By \thetag{\ref{LCpfaffian}} we have 
$\Pf({\widehat A}^{\overline I}_{\overline I})
=(-1)^{|\overline I|-M'+r'}\Pf(A)^{M'-r'-1}\Pf(A^I_I)$
and
$\Pf({\widehat B}^{\overline J}_{\overline J})
=(-1)^{|\overline J|-N'+l'+r'}\Pf(B)^{N'-l'-r'-1}\Pf(B^J_J)$.
Further it is easy to see that $s(I,\overline I)\equiv |I|-r'$ ($\MOD2$) 
and $s(J,\overline J)\equiv |J|-l'-r'$ ($\MOD2$)
Since $I\cup\overline I=R$ and $J\cup\overline J=S$,
we have $|I|+|\overline I|=\binom{M+1}{2}\equiv M'$ ($\MOD2$) 
and $|J|+|\overline J|=\binom{N+1}{2}\equiv N'$ ($\MOD2$).
These identities immediately imply \thetag{\ref{row_column_even}}.
$\Box$

%% file: lmsf6.tex
%
%
%
%
%
%


\section{Kawanaka's $q$-Littlewood formula}

In \cite{Ka1}, Kawanaka gave a certain $q$-series identity which is a generalization of the classical Schur-Littlewood identity
(see \cite{S}):
\begin{equation*}
\sum_{\lambda}s_{\lambda}(x)=\prod_{i}\frac1{1-x_i}\prod_{i<j}\frac1{1-x_ix_j}
\end{equation*}
where the sum runs over all partitions $\lambda$.
The Schur functions are well-known symmetric functions.
The reader should consult \cite{Ma} to see the detailed explanations of the symmetric functions.
Here we only use a well-known determinant expression for the Schur functions.
We use the notation in Macdonald's book \cite{Ma}.
For example, a \defterm{partition} is a non-increasing sequence of nonnegative integers $\lambda=(\lambda_1,\lambda_2,\dots)$ with finite non-zero parts.
The number of non-zero parts are called the \defterm{length} and denoted by $\ell(\lambda)$.
We assume the number of the variables is finite, 
say $n$, and $x=(x_1,\dots,x_n)$.
Then the \defterm{Schur function} corresponding to a partition $\lambda$ is defined to be
\begin{equation*}
s_{\lambda}(x)=\frac1{\Delta(x)}\det(x_{i}^{\lambda_j+n-j}).
\end{equation*}
Here $\Delta(x)=\prod_{i<j}(x_i-x_j)$.
In the following discussions we identify a partition with its Ferrers graph.
Given a partition $\lambda$,
the \defterm{hook-length} of $\lambda$ at $\alpha=(i,j)$ is, by definition, 
$h(\alpha)=h(i,j)=\lambda_i+\lambda_j-i-j+1$.
Let $(a;q)_{\infty}=\prod_{n=0}^{\infty}(1-aq^{n})$ and $(a;q)_{n}=\frac{(a;q)_{\infty}}{(aq^{n};q)_{\infty}}$ for complex numbers $a$, $q$ such that $|q|<1$.
We write $(a)_{n}$ (resp. $(a)_{\infty}$) for $(a;q)_{n}$ (resp. $(a;q)_{\infty}$) in short when there is no fear of confusion.

First of all, we recall Kawanaka's generalization of the Schur-Littlewood identity.
\begin{thm}(Kawanaka)
\label{kawanaka}
\begin{equation}
\label{kawa_id}
\sum_{\lambda}\prod_{\alpha\in\lambda}\frac{1+q^{h(\alpha)}}{1-q^{h(\alpha)}}s_{\lambda}(x)
=\prod_{i=1}^{n}\frac{(-x_{i}q;q)_{\infty}}{(x_{i};q)_{\infty}}\prod_{1\leq i<j\leq n}\frac{1}{1-x_ix_j}
\end{equation}
where the sum runs over all partitions $\lambda$.
\end{thm}

In this section we give a short proof of this identity as an 
application of the minor summation formula
and then use this method to obtain a similar formula as follows.
\begin{thm}
\label{ours}
\begin{equation}
\label{general_ours}
\sum_{\lambda}q^{n(\lambda)}\frac{\prod_{i=1}^{n}(a;q)_{\lambda_i+n-i}}{\prod_{\alpha\in\lambda}(1-q^{h(\alpha)})}s_{\lambda}(x)
=\prod_{i=1}^{n-1}(a;q)_{i}\prod_{i=1}^{n}\frac{(aq^{n-1}x_i;q)_{\infty}}{(x_{i};q)_{\infty}}.
\end{equation}
Here $n(\lambda)=\sum_{i\geq1}(i-1)\lambda_i$ for a partition $\lambda$.
\end{thm}
(The referee pointed out the second theorem is a consequence of Cauchy's identity and the specialization
of the Schur functions given in \cite{Ma} ch. I l.3 ex.3).
In order to prove the theorems, we first recall the $q$-binomial formula:
\begin{lem}
\label{q_binomial}
\begin{equation}
\label{qbinomial}
\sum_{n=0}^{\infty}\frac{(a)_{n}}{(q)_{n}}x^{n}=\frac{(ax)_{\infty}}{(x)_{\infty}}.
\Box
\end{equation}
\end{lem}
The following lemma is a generalization of the $q$-binomial formula and 
becomes the key to the proof of Kawanaka's identity.
\begin{lem}
\label{two_variables}
\begin{equation*}
\label{sum}
\sum_{k,l\geq0}
\frac{(-q)_{k}(-q)_{l}}{(q)_{k}(q)_{l}}
\frac{q^{l}-q^{k}}{q^{k}+q^{l}}x^{k}y^{l}
=\frac{(-qx)_{\infty}}{(x)_{\infty}}\frac{(-qy)_{\infty}}{(y)_{\infty}}
\frac{x-y}{1-xy}
\end{equation*}
\end{lem}
{\it Proof}.
Put
\begin{align*}
&F(x,y)=\sum_{k,l\geq0}\frac{(-q)_{k}(-q)_{l}}{(q)_{k}(q)_{l}}\frac{q^{l}-q^{k}}{q^{k}+q^{l}}x^{k}y^{l}=\sum_{k,l\geq0}a_{kl}x^ky^l,
\\
&G(x,y)=\frac{x-y}{1-xy}\prod_{r=0}^{\infty}\frac{1+xq^{r+1}}{1-xq^{r}}\frac{1+yq^{r+1}}{1-yq^{r}}=\sum_{k,l\geq0}b_{kl}x^ky^l.
\end{align*}
Then
\begin{equation*}
F(x,0)=\sum_{k\geq1}\frac{(-q)_{k-1}}{(q)_{k-1}}x^k=x\frac{(-qx)_{\infty}}{(x)_{\infty}}=G(x,0).
\end{equation*}
Exactly the same argument leads to $F(0,y)=-y\frac{(-qy)_{\infty}}{(y)_{\infty}}=G(0,y)$,
and This implies that $a_{k,0}=b_{k,0}$ and $a_{0,l}=b_{0,l}$ for $k,l\geq0$.
Next we claim that the coefficient of $x^ky^l$ of $(1-xy)F(x,y)$ is equal to the coefficient of $x^ky^l$ of $(1-xy)G(x,y)$ for $k,l\geq1$.
An easy calculation shows that
\begin{align*}
a_{kl}-a_{k-1,l-1}=2(q^l-q^k)\frac{(-q)_{k-1}(-q)_{l-1}}{(q)_{k}(q)_{l}}.
\end{align*}
On the other hand,
the coefficient of $x^ky^l$ in $(1-xy)G(x,y)=(x-y)\frac{(-qx)_{\infty}(-qy)_{\infty}}{(x)_{\infty}(y)_{\infty}}$ is 
\begin{align*}
\frac{(-q)_{k-1}}{(q)_{k-1}}\frac{(-q)_{l}}{(q)_{l}}-\frac{(-q)_{k}}{(q)_{k}}\frac{(-q)_{l-1}}{(q)_{l-1}}
=2(q^l-q^k)\frac{(-q)_{k-1}(-q)_{l-1}}{(q)_{k}(q)_{l}}.
\end{align*}
This shows our claim holds,
and in consequence we obtain $a_{kl}-a_{k-1,l-1}=b_{kl}-b_{k-1,l-1}$.
This proves the lemma by induction.
$\Box$
\bigbreak

A key observation to prove Kawanaka's formula is the following identity
which can be obtained from (1.7) of \cite{Ma}:
\begin{equation}
\prod_{\alpha\in\lambda}\frac{1+q^{h(\alpha)}}{1-q^{h(\alpha)}}
=\prod_{i=1}^{n}\frac{(-q)_{\ell_i}}{(q)_{\ell_i}}\prod_{i<j}\frac{1-q^{\ell_{i}-\ell_{j}}}{1+q^{\ell_{i}-\ell_{j}}}
.\label{key_id}
\end{equation}
Here $\ell_i=\lambda_i+n-i$ with $\ell(\lambda)\leq n$.
In fact,
in ex.~1, ch. I, l.1 of \cite{Ma},
the following identity is shown:
\begin{equation}
\prod_{\alpha\in\lambda}(1-q^{h(\alpha)})
=\frac
{\prod_{i=1}^{n}(q)_{\ell_i}}
{\prod_{i<j}\left(1-q^{\ell_{i}-\ell_{j}}\right)}.
\label{Macd_id}
\end{equation}
By the same argument one obtains
$$
\prod_{\alpha\in\lambda}(1+q^{h(\alpha)})
=\frac
{\prod_{i=1}^{n}(-q)_{\ell_i}}
{\prod_{i<j}\left(1+q^{\ell_{i}-\ell_{j}}\right)}
$$
from (1.7) of \cite{Ma}.
Taking the ratio of the equations one proves \thetag{\ref{key_id}}.
We also use the following famous identities.
(For the proof, see \cite{Ste}.)
\begin{lem}
Let $n$ be an even integer.
Let $x_1,\dots,x_n$ be indeterminates.
Then
\begin{align}
\label{PF1}
&\Pf\left[\frac{x_i-x_j}{x_i+x_j}\right]_{1\leq i<j\leq n}=\prod_{1\leq i<j\leq n}\frac{x_i-x_j}{x_i+x_j},\\
\label{PF2}
&\Pf\left[\frac{x_i-x_j}{1-x_ix_j}\right]_{1\leq i<j\leq n}=\prod_{1\leq i<j\leq n}\frac{x_i-x_j}{1-x_ix_j}.
\quad\Box
\end{align}
\end{lem}
\bigbreak

Let $A=(\alpha_{ij})_{i,j\geq0}$ denote a skew symmetric matrix.
As an application of Theorem~\ref{lmsf} 
we obtain the following formula from the definition of the Schur functions.
\begin{lem}
\label{Schur_sum}
Let $n$ be an even integer. 
We denote by 
$s_{\lambda}(x)$ the Schur functions of $n$ variables 
corresponding to a partition $\lambda$.
Then
\begin{equation}
\label{schur_func}
\sum_{\lambda}\Pf(\alpha_{\ell_p\ell_q})_{1\leq p,q\leq n}s_{\lambda}(x)=\frac{1}{\Delta(x)}\Pf(\beta_{ij})_{1\leq i,j\leq n}
\end{equation}
where $\lambda$ runs all the partition such that $\ell(\lambda)\leq n$,
\begin{equation*}
\beta_{ij}
=\sum_{k,l\geq0}\alpha_{kl}x_{i}^{k}x_{j}^{l}
=\sum_{0\leq k<l}\alpha_{kl}\begin{vmatrix}
x_{i}^{k}&x_{i}^{l}\\
x_{j}^{k}&x_{j}^{l}
\end{vmatrix},
\end{equation*}
and $\Delta(x)=\prod_{1\leq i<j\leq n}(x_{i}-x_{j})$.
\end{lem}
Now we are in position to prove Kawanaka's formula.
\bigbreak

\noindent
{\it Proof of Theorem~\ref{kawanaka}.}
It is enough to prove the case where $n$ is even.
For a partition $\lambda=(\lambda_1,\cdots,\lambda_n)$ ($\lambda_1\geq\cdots\geq\lambda_n\geq0$),
we put $\ell_i=\lambda_{i}+n-i$ ($i=1,2,\cdots,n$) as above.
If we put $\alpha_{kl}=\frac{(-q)_{k}(-q)_{l}}{(q)_{k}(q)_{l}}\frac{q^{l}-q^{k}}{q^{k}+q^{l}}$ in \thetag{\ref{schur_func}},
then \thetag{\ref{key_id}} and \thetag{\ref{PF1}} imply
\begin{align*}
\Pf\begin{bmatrix}
\alpha_{\ell_i\ell_j}
\end{bmatrix}_{1\leq i,j\leq n}
=\prod_{i=1}^{n}\frac{(-q)_{l_{i}}}{(q)_{l_{i}}}
\prod_{1\leq i<j\leq n}\frac{q^{l_{j}}-q^{l_{i}}}{q^{l_{i}}+q^{l_{j}}}
=\prod_{\alpha\in\lambda}\frac{1+q^{h(x)}}{1-q^{h(x)}}.
\end{align*}
Thus, by Lemma~\ref{two_variables} and \thetag{\ref{PF2}} we obtain
\begin{align*}
\Delta(x)\sum_{\lambda}\prod_{\alpha\in\lambda}\frac{1+q^{h(x)}}{1-q^{h(x)}}s_{\lambda}(x)
&=\Pf\begin{bmatrix}
\frac{(-qx_{i})_{\infty}}{(x_{i})_{\infty}}\frac{(-qx_{j})_{\infty}}{(x_{j})_{\infty}}
\frac{x_{i}-x_{j}}{1-x_{i}x_{j}}
\end{bmatrix}_{1\leq i<j\leq n}\\
&=\prod_{i=1}^{n}\frac{(-qx_{i})_{\infty}}{(x_{i})_{\infty}}
\Pf\begin{bmatrix}
\frac{x_{i}-x_{j}}{1-x_{i}x_{j}}
\end{bmatrix}_{1\leq i<j\leq n}\\
&=\Delta(x)\prod_{i=1}^{n}\frac{(-qx_{i})_{\infty}}{(x_{i})_{\infty}}
\prod_{1\leq i<j\leq n}
\frac{1}{1-x_{i}x_{j}}.
\end{align*}
This proves the theorem.
$\Box$

As a formula similar to \thetag{\ref{PF1}} and \thetag{\ref{PF2}},
the following is a Pfaffian version of the Vendermonde determinant.
\begin{lem}
\label{Vandermode_pfaffian}
Let $n=2r$ be an even integer.
Then
\begin{equation}
\label{vandermode_pf}
\Pf\left[\frac{(x_i^r-x_j^r)^2}{x_i-x_j}\right]_{1\leq i<j\leq n}
=\prod_{1\leq i<j\leq n}(x_i-x_j).\quad\Box
\end{equation}
\end{lem}

\bigbreak
This proof shows that replacing the entries of the Pfaffian by 
appropriate polynomials will be an interesting problem.
\bigbreak
\noindent
{\it Proof of Theorem~\ref{ours}}.
Now we consider a skew symmetric matrix $A=(\alpha_{k,l})$ of size $n=2r$ whose $(k,l)$-entry is defined by
\begin{equation*}
\alpha_{k,l}=\frac{(a)_{k}(a)_{l}}{(q)_{k}(q)_{l}}\frac{(q^{rl}-q^{rk})^2}{q^l-q^k}.
\end{equation*}
For a partition $\lambda=(\lambda_1,\dots,\lambda_n)$ we put $\ell_i=\lambda+n-i$. Then, from Lemma~\ref{Vandermode_pfaffian} we obtain
\begin{align*}
\Pf[\alpha_{\ell_i,\ell_j}]_{1\leq i<j\leq n}
&=\prod_{i=1}^{n}\frac{(a)_{i}}{(q)_{i}}\prod_{1\leq i<j\leq n}(q^{\ell_j}-q^{\ell_i})\\
&=q^{\sum_{i=1}^{n}(i-1)(\lambda_i+n-i)}\prod_{i=1}^{n}\frac{(a)_{i}}{(q)_{i}}\prod_{1\leq i<j\leq n}(1-q^{\ell_i-\ell_j})\\
&=\frac{q^{n(\lambda)+\frac16n(n-1)(n-2)}\prod_{i=1}^{n}(a)_{\lambda_i+n-i}}{\prod_{c\in\lambda}(1-q^{h(c)})}.
\end{align*}
Now, to apply Lemma~\ref{Schur_sum}, we need to study the sum:
\begin{align*}
f_{n}(x,y)
&=\sum_{k,l\geq0}\frac{(a)_{k}(a)_{l}}{(q)_{k}(q)_{l}}\frac{(q^{rl}-q^{rk})^2}{q^l-q^k}x^ky^l,\\
&=\sum_{k,l\geq0}\frac{(a)_{k}(a)_{l}}{(q)_{k}(q)_{l}}\left\{\sum_{\nu=1}^{r}q^{(\nu-1)k}q^{(n-\nu)l}-\sum_{\nu=1}^{r}q^{(n-\nu)k}q^{(\nu-1)l}\right\}x^ky^l,\\
&=\sum_{\nu=1}^{r}\frac{(aq^{\nu-1}x)_{\infty}}{(q^{\nu-1}x)_{\infty}}\frac{(aq^{n-\nu}y)_{\infty}}{(q^{n-\nu}y)_{\infty}}
-\sum_{\nu=1}^{r}\frac{(aq^{n-\nu}x)_{\infty}}{(q^{n-\nu}x)_{\infty}}\frac{(aq^{\nu-1}y)_{\infty}}{(q^{\nu-1}y)_{\infty}}.
\end{align*}
Thus we have
\begin{equation*}
f_{n}(x,y)=\frac{(aq^{n-1}x)_{\infty}(aq^{n-1}y)_{\infty}}{(x)_{\infty}(y)_{\infty}}g_{n}(x,y),
\end{equation*}
where $g_n(x,y)$ is a polynomial of $x$ and $y$ defined by
\begin{align*}
g_n(x,y)
&=\sum_{\nu=1}^{r}\prod_{k=1}^{\nu-1}(1-q^{k-1}x)\prod_{k=\nu}^{n-1}(1-aq^{k-1}x)\prod_{k=1}^{n-\nu}(1-q^{k-1}y)\prod_{k=n-\nu+1}^{n-1}(1-aq^{k-1}y)\\
&-\sum_{\nu=1}^{r}\prod_{k=1}^{n-\nu}(1-q^{k-1}x)\prod_{k=n-\nu+1}^{n-1}(1-aq^{k-1}x)\prod_{k=1}^{\nu-1}(1-q^{k-1}y)\prod_{k=\nu}^{n-1}(1-aq^{k-1}y).
\end{align*}
By applying Lemma~\ref{Schur_sum}, 
it is not hard to see that 
in order to prove \thetag{\ref{general_ours}},
it suffices to prove the following lemma, Lemma 6.8.
$\Box$

\begin{lem}
Let $n$ be even integer and let $g_n(x,y)$ be as above.
Then we have
\begin{equation*}
\Pf[g_n(x_i,x_j)]_{1\leq i<j\leq n}
=q^{\frac16n(n-1)(n-2)}\prod_{k=1}^{n-1}(a)_{k}\prod_{1\leq i<j\leq n}(x_i-x_j).
\end{equation*}
\end{lem}
{\it Proof}.
The method we use here is quite similar to that of we used in the proof of Lemma~\ref{Vandermode_pfaffian}.
First of all the reader should notice that $g(x,y)$ is of degree $(n-1)$ as a polynomial in the variable $x$.
This shows that $\Pf[g(x_i,x_j)]$ is a polynomial of degree at most $(n-1)$
if we see it as a polynomial of a fixed variable $x_i$.
Since $g(x,y)$ is skew symmetric,i.e. $g(y,x)=-g(x,y)$,
and this show that $(x_i-x_j)$ divides the Pfaffian,
and, as before, the complete produce $\prod_{i<j}(x_i-x_j)$ must divide
 the Pfaffian.
Thus we conclude that 
\begin{equation*}
\Pf[g(x_i,x_j)]_{1\leq i<j\leq n}
=c\prod_{1\leq i<j\leq n}(x_i-x_j).
\end{equation*}
If we see the right-hand side as a polynomial of a fixed variable $x_i$,
then it is of degree $(n-1)$, and this shows the constant $c$ must not include $x_i$.
Now, to determine the constant $c$,
which is independent of $x_i$,
we compare the coefficient of the monomial 
$\prod_{i=1}^{n}x_i^{n-i}$ of the both sides.
First we consider the left-hand side.
The Pfaffian is the sum of polynomials $\sgn\sigma g_n(x_{i_1},x_{j_1})\cdots g_n(x_{i_r},x_{j_r})$
for all perfect matching $\sigma=((i_1,j_1),\dots,(i_r,j_r))$ of $[n]$,
The monomial which contributes to the monomial $\prod_{i=1}^{n}x_i^{n-i}$ 
in the polynomial $g_n(x_{i_k},x_{j_k})$ is $x_{i_k}^{n-i_k}x_{j_k}^{n-j_k}$.
This shows hence that the coefficient of $\prod_{i=1}^{n}x_i^{n-i}$ in 
the left-hand side Pfaffian is equal to 
$\Pf\left[[x^{n-i}y^{n-j}]g_n(x,y)\right]$.
Thus, to prove the desired identity it suffices to prove the identity 
\begin{equation*}
\Pf\left[[x^{n-i}y^{n-j}]g_n(x,y)\right]_{1\leq i<j\leq n}=q^{\frac16n(n-1)(n-2)}\prod_{k=1}^{n-1}(a)_{k}.
\end{equation*}
Here $[x^ay^b]f(x,y)$ stands for the coefficient of 
the $x^ay^b$ in the polynomial $f(x,y)$.
Since the determination of the sign of the Pfaffain is easy, 
to prove the identity it suffices to show the following lemma.
$\Box$

\begin{lem}
Let $n=2r$ be an even integer and let $h_n(x,y)=h_n(a,b,q,t;x,y)$ be the polynomial of $x$ and $y$ defined by
\begin{align*}
&\sum_{\nu=1}^{r}\prod_{k=1}^{\nu-1}(1-q^{k-1}x)\prod_{k=\nu}^{n-1}(1-aq^{k-1}x)\prod_{k=1}^{n-\nu}(1-t^{k-1}y)\prod_{k=n-\nu+1}^{n-1}(1-bt^{k-1}y)\\
&-\sum_{\nu=1}^{r}\prod_{k=1}^{n-\nu}(1-q^{k-1}x)\prod_{k=n-\nu+1}^{n-1}(1-aq^{k-1}x)\prod_{k=1}^{\nu-1}(1-t^{k-1}y)\prod_{k=\nu}^{n-1}(1-bt^{k-1}y).
\end{align*}
Then we have
\begin{equation*}
\det\left[[x^{n-i}y^{n-j}]h_n(x,y)\right]_{1\leq i,j\leq n}=(qt)^{\frac16n(n-1)(n-2)}\prod_{i=1}^{n-1}(a;q)_i\prod_{j=1}^{n-1}(b;t)_j.
\end{equation*}
\end{lem}
{\it Proof}.
Note that $[x^{n-i}y^{n-j}]h_n(x,y)$ has degree $(n-1)$ regarded as a
polynomial in either $a$ or $b$,
This means the determinant has degree at most $n(n-1)$ in either variable $a$ or $b$.
We claim that $\prod_{i=1}^{n-1}(a;q)_i\prod_{j=1}^{n-1}(b;t)_j$ divides the determinant.
For this purpose we want to show that $(1-aq^{i-1})$ divides the determinant $(n-i)$ times for each $i=1,\dots,(n-1)$.
This can be done by substituting $a=q^{-k}$ ($k=0,\dots,n-2$) into the determinant and computing the rank.
The details are left to the reader.
By this argument, one see that
\begin{equation*}
\det\left[[x^{n-i}y^{n-j}]h_n(x,y)\right]_{1\leq i,j\leq n}=c\prod_{i=1}^{n-1}(a;q)_i\prod_{j=1}^{n-1}(b;t)_j,
\end{equation*}
where $c$ is a constant independent of $a$ and $b$.
To find $c$, Compare the constant term of the both sides regarding them as polynomials of $a$ and $b$.
$\Box$

%% file: lmsf7.tex
%
%
%
%
%
%


\section{Kawanaka's $q$-Cauchy identity}

In \cite{Ka2}
Kawanaka gave a $q$-Cauchy formula,
which is regarded as 
a determinant version of Kawanaka's $q$-Littlewood formula in the previous section.
Before we state the theorem we need some definitions.
Let $\lambda$ and $\mu$ be partitions,
and let $c=(i,j)$ be any cell in the plane.
As a natural generalization of the ordinary hook length $h_{\lambda}(c)=\lambda_i+\lambda_j'-i-j+1$,
Kawanaka introduced a new statistic
\begin{equation*}
h_{\lambda\mu}(c)=\lambda_i+\mu_j'-i-j+1
\end{equation*}
in \cite{Ka2}.
For example, let $\lambda=(4,3,1,1)$ and $\mu=(3,3)$.
If we fill each cell $c$ of $\lambda$ with the numbers $h_{\lambda\mu}(c)$,
then it looks as follows:
$$
\young{
5&4&3&0\\
3&2&1\\
0\\
-1\\
}
$$
In \cite{Ka2}, he also defined
\begin{equation*}
n(\lambda,\mu)
=\sum_{(i,j)\in\lambda-\mu}(\lambda_{j}'-i)
=\sum_{(i,j)\in\lambda-\mu}(i-\mu_{j}'-1)
\end{equation*}
which is regarded as a generalization of $n(\lambda)=\sum_{i\geq1}(i-1)\lambda_{i}$ in \cite{Ma}.
Let $t$ be an indeterminate.
For any partitions $\lambda$ and $\mu$,
define a rational function $J_{\lambda\mu}(t)$ in $t$ by
\begin{eqnarray*}
&&J_{\lambda\mu}(t)
=t^{n(\lambda,\mu)}\prod_{c\in\lambda}\frac{1+t^{h_{\lambda\mu}(c)}}{1-t^{h_{\lambda}(c)}}
\cdot t^{n(\mu,\lambda)}\prod_{c\in\mu}\frac{1+t^{h_{\mu\lambda}(c)}}{1-t^{h_{\mu}(c)}}.
\end{eqnarray*}

\begin{thm} (Kawanaka)
\label{thm:det-kawanaka}
Let $x=(x_1,x_2,\dots)$ and $y=(y_1,y_2,\dots)$ be two independent sequences of variables.
For a partition $\lambda$ let $s_{\lambda}(x)$ and $s_{\lambda}(y)$ be the corresponding Schur functions in $x$ and $y$ respectively.
Then we have the following identity:
\begin{eqnarray}
&&\sum_{\lambda,\mu}q^{|\lambda-\mu|+|\mu-\lambda|}
J_{\lambda\mu}(q^2)s_{\lambda}(x)s_{\mu}(y)
\nonumber\\
&&=
\prod_{i\geq 1}\frac{(-qx_i;q^2)_{\infty}}{(qx_i;q^2)_{\infty}}
\prod_{j\geq 1}\frac{(-qy_j;q^2)_{\infty}}{(qy_j;q^2)_{\infty}}
\prod_{i,j\geq 1}\frac1{1-x_iy_j}.
\label{eqn:kawanaka}
\end{eqnarray}
Here $|\lambda-\mu|$ is the number of the cells in the set-theoretical difference $\{c:c\in\lambda,c\not\in\mu\}$
and $|\mu-\lambda|$ is the number of the cells in $\{c:c\not\in\lambda,c\in\mu\}$.
\end{thm}

%
%
\begin{lem}
\label{lemma:two-vars}
\begin{align}
\sum_{k,l\geq0}
\frac{(-q^2;q^2)_{k}}{(q^2;q^2)_{k}}\frac{(-q^2;q^2)_{l}}{(q^2;q^2)_{l}}
\frac{2x^ky^l}{q^{k-l}+q^{l-k}}
=\frac{(-qx;q^2)_{\infty}}{(qx;q^2)_{\infty}}\frac{(-qy;q^2)_{\infty}}{(qy;q^2)_{\infty}}
\frac1{1-xy}.
\end{align}
\end{lem}
%
%
{\it Proof}.
Put
\begin{align*}
&F(x,y)=\sum_{k,l\geq0}\frac{(-q^2;q^2)_{k}}{(q^2;q^2)_{k}}\frac{(-q^2;q^2)_{l}}{(q^2;q^2)_{l}}
\frac{2x^ky^l}{q^{k-l}+q^{l-k}}
=\sum_{k,l\geq0}a_{kl}x^ky^l,
\\
&G(x,y)
=\frac{(-qx;q^2)_{\infty}}{(qx;q^2)_{\infty}}
\frac{(-qy;q^2)_{\infty}}{(qy;q^2)_{\infty}}\frac1{1-xy}
.
\end{align*}
First, we compare the coefficients of $x^ky^l$ in $(1-xy)F(x,y)$ and $(1-xy)G(x,y)$.
By Lemma~\ref{q_binomial} we have
\begin{equation*}
(1-xy)G(x,y)=\sum_{k,l\geq0}\frac{(-1;q^2)_{k}(-1;q^2)_{l}}{(q^2;q^2)_{k}(q^2;q^2)_{l}}q^{k+l}x^ky^l.
\end{equation*}
Meanwhile, the coefficient of $x^ky^l$ in $(1-xy)F(x,y)$ for $k,l\geq1$ is equal to
\begin{align*}
&a_{kl}-a_{k-1,l-1}\\
&=\frac{(-q^2;q^2)_{k-1}(-q^2;q^2)_{l-1}}{(q^2;q^2)_{k}(q^2;q^2)_{l}}
\{(1+q^{2k})(1+q^{2l})-(1-q^{2k})(1-q^{2l})\}\frac{2q^kq^l}{q^{2k}+q^{2l}}\\
&=4\frac{(-q^2;q^2)_{k-1}(-q^2;q^2)_{l-1}}{(q^2;q^2)_{k}(q^2;q^2)_{l}}q^{k+l}
=\frac{(-1;q^2)_{k}(-1;q^2)_{l}}{(q^2;q^2)_{k}(q^2;q^2)_{l}}q^{k+l}.
\end{align*}
When $l=0$,
it is easy to see that the coefficient of $x^k$ in $(1-xy)F(x,y)$ is equal to
\begin{equation*}
\frac{(-q^2;q^2)_{k}}{(q^2;q^2)_{k}}\frac{2q^k}{1+q^{2k}}=\frac{(-1;q^2)_{k}}{(q^2;q^2)_{k}}q^k.
\end{equation*}
On the other hand,
when $=0$,
it is also easy to see that the coefficient of $y^l$ in $(1-xy)F(x,y)$ is equal to $\frac{(-1;q^2)_{l}}{(q^2;q^2)_{l}}q^l$.
Thus, 
the coefficients agree in all cases,
and we conclude that $(1-xy)F(x,y)=(1-xy)G(x,y)$.
This completes the proof.
$\Box$
\medbreak

The following identities are known as the Cauchy determinants (see \cite{Ma}).
%
%
\begin{prop}
\label{Cauchy_det}
\begin{align}
\label{Cauchy1}
&\det\left(\frac1{x_i+y_j}\right)_{1\leq i,j\leq n}
=\frac{\prod_{1\leq i<j\leq n}(x_i-x_j)\prod_{1\leq i<j\leq n}(y_i-y_j)}{\prod_{1\leq i,j\leq n}(x_i+y_j)},
\\
\label{Cauchy2}
&\det\left(\frac1{1-x_iy_j}\right)_{1\leq i,j\leq n}
=\frac{\prod_{1\leq i<j\leq n}(x_i-x_j)\prod_{1\leq i<j\leq n}(y_i-y_j)}{\prod_{1\leq i,j\leq n}(1-x_iy_j)}.
\end{align}
\end{prop}
%
%
%
%
\begin{lem}
\label{lemma:coeff}
Let $\lambda=(\lambda_1,\dots,\lambda_n)$ and $\mu=(\mu_1,\dots,\mu_n)$ be partitions such that $\ell(\lambda),\ell(\mu)\leq n$.
We put $k_i=\lambda_{i}+n-i$ and $\ell_i=\mu_{i}+n-i$ for $1\leq i\leq n$.
If we put
\begin{align*}
a_{kl}=\frac{(-q^2;q^2)_{k}}{(q^2;q^2)_{k}}\frac{(-q^2;q^2)_{l}}{(q^2;q^2)_{l}}\frac{2q^{k+l}}{q^{2k}+q^{2l}},
\end{align*}
then
\begin{equation}
\det(a_{k_i\ell_j})_{1\leq i,j\leq n}
=q^{|\lambda-\mu|+|\mu-\lambda|}J_{\lambda\mu}(q^2).
\end{equation}
\end{lem}
%
%
{\it Proof}.
From \thetag{\ref{Cauchy1}} we obtain 
\begin{align*}
&\det(a_{k_i\ell_j})_{1\leq i,j\leq n}\\
&=2^nq^{\sum_{i}(2i-1)k_i+\sum_{j}(2j-1)\ell_j}
\frac{\prod_{i<j}(1-q^{2(k_i-k_j)})}{\prod_{i}(q^2;q^2)_{k_i}}\cdot
\\
&\qquad\times\frac{\prod_{i<j}(1-q^{2(\ell_i-\ell_j)})}{\prod_{j}(q^2;q^2)_{\ell_j}}\cdot
\frac{\prod_{i}(-q^2;q^2)_{k_i}\prod_{j}(-q^2;q^2)_{\ell_j}}{\prod_{i,j}(q^{2k_i}+q^{2\ell_j})}
.
\end{align*}
By \thetag{\ref{Macd_id}} we have
\begin{align*}
\prod_{c\in\lambda}\frac1{1-q^{2h_{\lambda}(c)}}=\frac{\prod_{1\leq i<j\leq n}(1-q^{2(k_i-k_j)})}{\prod_{i=1}^n (q^2;q^2)_{k_i}},
\\
\prod_{c\in\lambda}\frac1{1-q^{2h_{\mu}(c)}}=\frac{\prod_{1\leq i<j\leq n}(1-q^{2(\ell_i-\ell_j)})}{\prod_{j=1}^n (q^2;q^2)_{\ell_j}}.
\end{align*}
%
Thus it is enough to show that
\begin{eqnarray*}
&&2^nq^{\sum_{i=1}^{n}(2i-1)k_i+\sum_{j=1}^{n}(2j-1)\ell_j}
\frac{\prod_{i=1}^{n}(-q^2;q^2)_{k_i}\prod_{j=1}^{n}(-q^2;q^2)_{\ell_j}}{\prod_{i,j=1}^{n}(q^{2k_i}+q^{2\ell_j})}
\\
&&=q^{|\lambda-\mu|+|\mu-\lambda|+2n(\lambda,\mu)+2n(\mu,\lambda)}
\prod_{c\in\lambda}(1+q^{2h_{\lambda\mu}(c)})
\prod_{c\in\mu}(1+q^{2h_{\mu\lambda}(c)}).
\end{eqnarray*}
Note that $h_{\lambda\mu}(i,j)=\lambda_i-j+\mu_j'-i+1=\lambda_{i}-j+\sharp\{r:\mu_r\geq j\}-i+1$.
Here $\sharp\,A$ stands for the cardinality of the set $A$.
For a fixed $i$,
let $n_{i}=\sharp\{r:\mu_r>\lambda_i\}$ denote the number of parts of $\mu$
which is greater than $\lambda_i$.
Thus we have $\sharp\{r:\mu_r\leq\lambda_i\}=n-n_{i}$.
Then,
if we draw the Young diagram of $\lambda$ and 
fill each cell $c$ with $h_{\lambda\mu}(c)$,
then we find that the numbers in the $i$th row of $\lambda$ are
\begin{equation*}
\left[n_{i}-i+1,k_i\right]-\{k_i-\ell_r:n_i<r\leq n\}.
\end{equation*}
Here we write $[a,b]=\{a,a+1,a+2,\dots,b\}$ for integers $a,b\in\Bbb{Z}$.
In fact fix an $i$, and put 
$\mu^{(i)}=(\mu^{(i)}_{1},\dots,\mu^{(i)}_{n-n_i})=(\mu_{n_{i}+1},\mu_{n_{i}+2},\dots,\mu_{n})$.
Note that the conjugate ${\mu^{(i)}}'$ fits in the box $\lambda_i\times(n-n_i)$.
Use \cite{Ma}, I. (1.7) to the partition ${\mu^{(i)}}'$,
then we obtain
$
\{{\mu^{(i)}_{r}}'+\lambda_i-r:1\leq r\leq\lambda_i\}
\uplus\{\lambda_i-1+r-\mu^{(i)}_{r}:1\leq r\leq n-n_i\}
=[0,\lambda_i+n-n_i-1].
$
Since ${\mu^{(i)}_{r}}'=\mu_r'-n_i$,
we conclude that
$
\{h_{\lambda\mu}(i,r):1\leq r\leq\lambda_i\}
\uplus\{k_i-\mu_{r+n_i}:1\leq r\leq n-n_i\}
=[n_i-i+1,k_i].
$

\begin{figure}[b]
\begin{center}
\setlength{\unitlength}{0.20mm}
\begin{picture}(180,160)
\thicklines
\put( 30,180){\line(1,0){30}}
\put( 60,180){\line(1,0){30}}
\put( 90,180){\line(1,0){30}}
\put(120,180){\line(1,0){30}}
\put(150,180){\line(1,0){30}}
\put(180,180){\line(1,0){30}}
\put(210,180){\line(1,0){30}}
\put( 30,150){\line(1,0){30}}
\put( 60,150){\line(1,0){30}}
\put( 90,150){\line(1,0){30}}
\put(120,150){\line(1,0){30}}
\put(150,150){\line(1,0){30}}
\put(180,150){\line(1,0){30}}
\put( 30,120){\line(1,0){30}}
\put( 60,120){\line(1,0){30}}
\put( 90,120){\line(1,0){30}}
\put(120,120){\line(1,0){30}}
\put( 30, 90){\line(1,0){30}}
\put( 60, 90){\line(1,0){30}}
\put( 90, 90){\line(1,0){30}}
\put(120, 90){\line(1,0){30}}
\put( 30, 60){\line(1,0){30}}
\put( 30, 30){\line(1,0){30}}
\put( 30,  0){\line(1,0){30}}
\put( 30,150){\line(0,1){30}}
\put( 60,150){\line(0,1){30}}
\put( 90,150){\line(0,1){30}}
\put(120,150){\line(0,1){30}}
\put(150,150){\line(0,1){30}}
\put(180,150){\line(0,1){30}}
\put(210,150){\line(0,1){30}}
\put( 30,120){\line(0,1){30}}
\put( 60,120){\line(0,1){30}}
\put( 90,120){\line(0,1){30}}
\put(120,120){\line(0,1){30}}
\put(150,120){\line(0,1){30}}
\put( 30, 90){\line(0,1){30}}
\put( 60, 90){\line(0,1){30}}
\put( 90, 90){\line(0,1){30}}
\put(120, 90){\line(0,1){30}}
\put(150, 90){\line(0,1){30}}
\put( 30, 60){\line(0,1){30}}
\put( 60, 60){\line(0,1){30}}
\put( 30, 30){\line(0,1){30}}
\put( 60, 30){\line(0,1){30}}
\put( 30,  0){\line(0,1){30}}
\put( 60,  0){\line(0,1){30}}
\put( 30,180){\circle*{10}}
\put( 60,180){\circle*{10}}
\put( 90,180){\circle*{10}}
\put(120,180){\circle*{10}}
\put(150,180){\circle*{10}}
\put(180,180){\circle*{10}}
\put(210,180){\circle*{10}}
\put(240,180){\circle*{10}}
\put( 30,150){\circle*{10}}
\put( 60,150){\circle*{10}}
\put( 90,150){\circle*{10}}
\put(120,150){\circle*{10}}
\put(150,150){\circle*{10}}
\put(180,150){\circle*{10}}
\put(210,150){\circle*{10}}
\put( 90,120){\circle*{10}}
\put(120,120){\circle*{10}}
\put(150,120){\circle*{10}}
\put( 90, 90){\circle*{10}}
\put(120, 90){\circle*{10}}
\put(150, 90){\circle*{10}}
\put( 30,120){\circle*{7}}
\put( 60,120){\circle*{7}}
\put( 30, 90){\circle*{7}}
\put( 60, 90){\circle*{7}}
\put( 30, 60){\circle*{7}}
\put( 60, 60){\circle*{7}}
\put( 30, 30){\circle*{7}}
\put( 60, 30){\circle*{7}}
\put( 30,  0){\circle*{7}}
\put( 60,  0){\circle*{7}}
\put( 35,167){$\scriptstyle\bold0$}
\put( 65,167){$\scriptstyle\bold1$}
\put( 95,167){$\scriptstyle\bold2$}
\put(125,167){$\scriptstyle\bold3$}
\put(155,167){$\scriptstyle\bold4$}
\put(185,167){$\scriptstyle\bold5$}
\put(215,167){$\scriptstyle\bold6$}
\put(245,167){$\scriptstyle\bold7$}
\put( 35,137){$\scriptstyle-\bold1$}
\put( 65,137){$\scriptstyle\bold0$}
\put( 95,137){$\scriptstyle\bold1$}
\put(125,137){$\scriptstyle\bold2$}
\put(155,137){$\scriptstyle\bold3$}
\put(185,137){$\scriptstyle\bold4$}
\put(215,137){$\scriptstyle\bold5$}
\put( 95,107){$\scriptstyle\bold0$}
\put(125,107){$\scriptstyle\bold1$}
\put(155,107){$\scriptstyle\bold2$}
\put( 95, 77){$\scriptstyle-\bold1$}
\put(125, 77){$\scriptstyle\bold0$}
\put(155, 77){$\scriptstyle\bold1$}
\put( 35,107){$\scriptstyle\it2$}
\put( 65,107){$\scriptstyle\it1$}
\put( 35, 77){$\scriptstyle\it3$}
\put( 65, 77){$\scriptstyle\it2$}
\put( 35, 47){$\scriptstyle\it4$}
\put( 65, 47){$\scriptstyle\it3$}
\put( 35, 17){$\scriptstyle\it5$}
\put( 65, 17){$\scriptstyle\it4$}
\put( 35,-13){$\scriptstyle\it6$}
\put( 65,-13){$\scriptstyle\it5$}
\end{picture}
\caption{Lattice points}
\label{figure:lattice}
\end{center}
\end{figure}
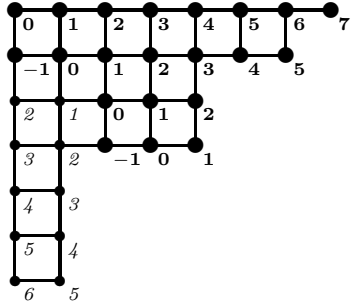

Put $m_j=\sharp\{r:\lambda_r\geq\mu_j\}$ for each $1\leq j\leq n$.
The same argument shows that,
if we write the Young diagram of $\mu$ and fill each cell $c$ 
with $h_{\mu\lambda}(c)$,
then the numbers in the $j$th row of $\mu$ are
\begin{equation*}
\left[m_j-j+1,\ell_j\right]-\{\ell_j-k_r:m_j<r\leq n\}.
\end{equation*}
Thus it is enough to show that
\begin{align*}
&2^nq^{\sum_{i=1}^{n}(2i-1)k_i+\sum_{j=1}^{n}(2j-1)\ell_j}
\frac{\prod_{i=1}^{n}(-q^2;q^2)_{k_i}\prod_{j=1}^{n}(-q^2;q^2)_{\ell_j}}{\prod_{i,j=1}^{n}(q^{2k_i}+q^{2\ell_j})}
\\
&=q^{|\lambda-\mu|+|\mu-\lambda|+2n(\lambda,\mu)+2n(\mu,\lambda)}
\frac{\prod_{i=1}^n\prod_{r=n_i-i+1}^{k_i}\left(1+q^{2r}\right)
\prod_{j=1}^n\prod_{r=m_j-j+1}^{\ell_j}\left(1+q^{2r}\right)}
{\prod_{\lambda_i\geq\mu_j}\left(1+q^{2(k_i-\ell_j)}\right)
\prod_{\lambda_i<\mu_j}\left(1+q^{2(\ell_i-k_j)}\right)}
\end{align*}
First it is easy to see that
\begin{eqnarray*}
\prod_{\mu_j\leq\lambda_i}(1+q^{2(k_i-\ell_j)})
\prod_{\lambda_i<\mu_j}(1+q^{2(\ell_i-k_j)})
=q^{-2P(\lambda,\mu)}\prod_{i,j=1}^{n}(q^{2k_i}+q^{2\ell_j}),
\end{eqnarray*}
where
\begin{eqnarray*}
P(\lambda,\mu)=\sum_{\lambda_i\geq\mu_j}\ell_j+\sum_{\lambda_i<\mu_j}k_i
=\sum_{i=1}^{n}n_ik_i+\sum_{j=1}^{n}m_j\ell_{j}.
\end{eqnarray*}
Next we claim that
\begin{align*}
\prod_{i=1}^n\prod_{r=n_i-i+1}^{k_i}\left(1+q^{2r}\right)
\prod_{j=1}^n\prod_{r=m_j-j+1}^{\ell_j}\left(1+q^{2r}\right)
&=q^{-2Q(\lambda,\mu)}2^n\prod_{i=1}^{n}(-q^2;q^2)_{k_i}\prod_{j=1}^{n}(-q^2;q^2)_{\ell_j},
\end{align*}
where
\begin{align*}
Q(\lambda,\mu)
&=
\sum_{{i}\atop{\lambda_i\geq\mu_i}}\binom{i-1-n_i}2
+\sum_{{i}\atop{\lambda_i<\mu_i}}\binom{i-1-m_i}2.
\end{align*}
In fact, let $A$ and $B$ be the sets of lattice points defined by
\begin{eqnarray*}
&&A=\bigcup_{i=1}^{n}\{(i-1,y):n_i\leq y\leq \lambda_i+n-1\},\\
&&B=\bigcup_{j=1}^{n}\{(x,j-1):m_j\leq x\leq \mu_j+n-1\}.
\end{eqnarray*}
Then we have
$\prod_{i=1}^n\prod_{r=n_i-i+1}^{k_i}\left(1+q^{2r}\right)
=\prod_{(x,y)\in A}(1+q^{2(y-x)})
$
and
$\prod_{j=1}^n\prod_{r=m_j-j+1}^{\ell_j}\left(1+q^{2r}\right)
=\prod_{(x,y)\in B}(1+q^{2(x-y)}).$
For example,
if $n=4$, $\lambda=(4,3,1,1)$ and $\mu=(3,3)$,
then the big circles in Figure~\ref{figure:lattice} are in $A$ and the small circles are in $B$.
The numbers assigned to big circles are $y-x$ and the numbers assigned to small circles are $x-y$.
Put
$A_{1}=\bigcup_{i=1}^{n}\{(i-1,y):n_i\leq y\leq n-1\}$,
$B_{1}=\bigcup_{j=1}^{n}\{(x,j-1):m_i\leq x\leq n-1\}$,
$A_{2}=\bigcup_{i=1}^{n}\{(i-1,y):n\leq y\leq \lambda_i+n-1\}$
and
$B_{2}=\bigcup_{j=1}^{n}\{(x,j-1):n\leq x\leq \mu_j+n-1\}$.
Then we have $A=A_1\cup A_2$ and $B=B_1\cup B_2$.
It is also easy to see that $A_1\cup B_1=[0,n-1]\times[0,n-1]$,
which implies that,
as a multi-set,
\begin{equation*}
\bigcup_{i=1}^{n}\{|y-i+1|:n_i\leq y\leq n-1\}
\cup
\bigcup_{j=1}^{n}\{|x-j+1|:m_i\leq x\leq n-1\}
\end{equation*}
is equal to $\{|x-y|;(x,y)\in[0,n-1]\times[0,n-1]\}$,
and is 
composed of $n$ $0$'s, $2(n-1)$ $1$'s, $2(n-2)$ $2$'s, \dots, $2$ $(n-1)$'s.
This shows that
\begin{eqnarray*}
&&\prod_{(x,y)\in A}(1+q^{2(y-x)})
\prod_{(x,y)\in B}(1+q^{2(x-y)})\\
&&=q^{-2Q(\lambda,\mu)}2^n\prod_{i=1}^{n}(-q^2;q^2)_{k_i}\prod_{j=1}^{n}(-q^2;q^2)_{\ell_j}.
\end{eqnarray*}
In another word we can restate
\begin{align*}
Q(\lambda,\mu)
&=
\sum_{{(x,y)\in A}\atop{x<y}}(y-x)+\sum_{{(x,y)\in B}\atop{x>y}}(x-y).
\end{align*}
Thus the proof will be done if we prove the following identity:
\begin{eqnarray*}
&&|\lambda-\mu|+|\mu-\lambda|+2n(\lambda,\mu)+2n(\mu,\lambda)+2P(\lambda,\mu)-2Q(\lambda,\mu)
\\
&&=\sum_{i=1}^n(2i-1)k_i+\sum_{j=1}^n(2j-1)\ell_j.
\end{eqnarray*}
In the above example,
we have $|\lambda-\mu|=3$, $|\mu-\lambda|=n(\mu,\lambda)=0$, $n(\lambda,\mu)=1$, $P(\lambda,\mu)=64$, $Q(\lambda,\mu)=4$,
and $3+2+64-4=65=\sum_{i=1}^4(2i-1)k_i+\sum_{j=1}^4(2j-1)\ell_j$.
In the following lemma we prove this identity.
$\Box$
%
%
\begin{lem}
Let $n$ be a nonnegative integer,
and let $\lambda$ and $\mu$ be partitions which satisfies $\ell(\lambda),\ell(\mu)\leq n$.
Let $P(\lambda,\mu)$ and $Q(\lambda,\mu)$ be as above.
Then the identity 
\begin{align}
&|\lambda-\mu|+|\mu-\lambda|+2n(\lambda,\mu)+2n(\mu,\lambda)+2P(\lambda,\mu)-2Q(\lambda,\mu)
\nonumber\\
&=\sum_{i=1}^n(2i-1)k_i+\sum_{j=1}^n(2j-1)\ell_j.
\label{id:strange}
\end{align}
holds.
\end{lem}
{\it Proof.}
We proceed by induction on $n$.
When $n=1$, assume $\lambda_1\geq\mu_1$.
Then it is easy to see that $n(\lambda,\mu)=n(\mu,\lambda)=Q(\lambda,\mu)=0$
and $P(\lambda,\mu)=\mu$.
This shows that the left-hand side equals $\lambda_1+\mu_1$
and it coincides with the right-hand sides.
In the case $\lambda_1<\mu_1$, we can prove it similarly.
Assume $n\geq2$ and \thetag{\ref{id:strange}} holds up to $(n-1)$.
Given partitions $\lambda=(\lambda_1,\dots,\lambda_{n})$ and $\mu=(\mu_1,\dots,\mu_{n})$,
we put $\widetilde\lambda=(\widetilde\lambda_1,\dots,\widetilde\lambda_{n-1})=(\lambda_1,\dots,\lambda_{n-1})$ and $\widetilde\mu=(\widetilde\mu_1,\dots,\widetilde\mu_{n-1})=(\mu_1,\dots,\mu_{n-1})$.
Further we set $\widetilde k_{i}=\widetilde\lambda_{i}+n-1-i$ and $\widetilde\ell_{i}=\widetilde\mu_{i}+n-1-i$ for $1\leq i\leq n-1$.
Then,
by the induction hypothesis,
we may assume that \thetag{\ref{id:strange}} holds for $(n-1)$, $\widetilde\lambda$, $\widetilde\mu$, $\widetilde k$ and $\widetilde\ell$.
First we assume that $\lambda_{n}\geq\mu_{n}$.
Thus we have $|\lambda-\mu|+|\mu-\lambda|=|\widetilde\lambda-\widetilde\mu|+|\widetilde\mu-\widetilde\lambda|+\lambda_{n}-\mu_{n}$.
From the condition $\lambda_n\geq\mu_n$,
we can find an integer $s$ such that $0\leq s<n$ and $\mu_{s}>\lambda_{n}\geq\mu_{s+1}$ holds.
Here we use the convention that $\lambda_0=\mu_0=\infty$.
Using this $s$,
we can express the statistics on $\lambda$ and $\mu$ with the statistics on $\widetilde\lambda$ and $\widetilde\mu$.
For example,
if we write the Young diagram of $\lambda$ and $\mu$
and fill the cell $(i,j)\in\lambda-\mu$ with the number $\mu_j'-i-1$,
then we easily see that
\begin{align*}
n(\lambda,\mu)&=n(\widetilde\lambda,\widetilde\mu)
+(n-s-1)(\lambda_{n}-\mu_{s+1})+\sum_{i=s+1}^{n-1}(n-i-1)(\mu_{i}-\mu_{i+1}),\\
&=n(\widetilde\lambda,\widetilde\mu)+(n-s-1)\lambda_{n}-\sum_{i=s+1}^{n-1}\mu_i\\
n(\mu,\lambda)&=n(\widetilde\mu,\widetilde\lambda).
\end{align*}
For a fixed $i$ such that $1\leq i<n$,
from the fact that $\mu_s>\lambda_n\geq\mu_{s+1}$,
it is easy to see that
\begin{align*}
&\sharp\{r:\mu_r>\lambda_i\}=\sharp\{r:\widetilde\mu_r>\widetilde\lambda_i\},\\
&\sharp\{r:\lambda_r\geq\mu_i\}
=\begin{cases}
\sharp\{r:\widetilde\lambda_r\geq\widetilde\mu_i\}&\text{ if $1\leq i\leq s$,}\\
\sharp\{r:\widetilde\lambda_r\geq\widetilde\mu_i\}+1&\text{ if $s+1\leq i< n$.}
\end{cases}
\end{align*}
From these facts we have 
\begin{align*}
&P(\lambda,\mu)=P(\widetilde\lambda,\widetilde\mu)
+(n-1)^2
+\sum_{j=s+1}^{n-1}\widetilde\ell_{j}+s\lambda_{n}+n\mu_{n},\\
&Q(\lambda,\mu)=Q(\widetilde\lambda,\widetilde\mu)
+\binom{n-1-s}2.
\end{align*}
Here we used the fact 
$\sum_{i=1}^{n-1}\sharp\{r:\mu_r>\lambda_i\}+\sum_{j=1}^{n-1}\sharp\{r:\lambda_r\geq\mu_j\}=(n-1)^2$,
which is easy to confirm. From these identities,
we obtain
\begin{align*}
&|\lambda-\mu|+|\mu-\lambda|+2n(\lambda,\mu)+2n(\mu,\lambda)+2P(\lambda,\mu)-2Q(\lambda,\mu)\\
&=|\widetilde\lambda-\widetilde\mu|+|\widetilde\mu-\widetilde\lambda|
+2n(\widetilde\lambda,\widetilde\mu)+2n(\widetilde\mu,\widetilde\lambda)+2P(\widetilde\lambda,\widetilde\mu)-2Q(\widetilde\lambda,\widetilde\mu)\\
&+(2n-1)\lambda_n+(2n-1)\mu_n
+2(n-1)^2.
\end{align*}
By the induction hypothesis we have 
$|\widetilde\lambda-\widetilde\mu|+|\widetilde\mu-\widetilde\lambda|
+2n(\widetilde\lambda,\widetilde\mu)+2n(\widetilde\mu,\widetilde\lambda)+2P(\widetilde\lambda,\widetilde\mu)-2Q(\widetilde\lambda,\widetilde\mu)
=\sum_{i=1}^{n-1}(2i-1)\widetilde k_{i}+\sum_{i=1}^{n-1}(2i-1)\widetilde \ell_{i}
=\sum_{i=1}^{n-1}(2i-1)k_{i}+\sum_{i=1}^{n-1}(2i-1)\ell_{i}-2(n-1)^2
$,
and this proves the desired identity.
In the case of $\lambda_n<\mu_n$,
we may find an integer $s$ which satisfies $0\leq s<n$ and $\lambda_{s}\geq\mu_n>\lambda_{s+1}$.
A similar argument will lead to the desired identity again.
$\Box$

\bigbreak
\noindent
{\it Proof of Theorem~\ref{thm:det-kawanaka}.}
We may assume that the number of variables are finite,
i.e., $x=(x_1,\dots,x_n)$ and $y=(y_1,\dots,y_n)$.
Assume $N\geq n$ is a positive integer.
Let $T$ and $S$ be two $n$ by $N$ rectangular matrices defined by
\begin{eqnarray*}
T=\left(x_i^{N-j}\right)^{i=1,\dots,n}_{j=1,\dots,N},
\qquad\qquad
S=\left(y_i^{N-j}\right)^{i=1,\dots,n}_{j=1,\dots,N}.
\end{eqnarray*}
Let $A$ be an $N$ by $N$ square matrix defined by
\begin{eqnarray*}
A=\left(\frac{(-q^2;q^2)_{N-i}}{(q^2;q^2)_{N-i}}\frac{(-q^2;q^2)_{N-j}}{(q^2;q^2)_{N-j}}\frac2{q^{i-j}+q^{j-i}}\right)_{i,j=1,\dots,N}.
\end{eqnarray*}
Now we compute $\lim_{N\rightarrow\infty}\det{}^tTAS$ in two different ways.
By the Cauchy-Binet formula \thetag{\ref{eq_gene_cauchy-binet}},
we have
\begin{equation*}
\det{}^tTAS
=\sum_{{I\subseteq[N]}\atop{\sharp I=n}}\sum_{{J\subseteq[N]}\atop{\sharp J=n}}\det T_I\det A^I_J \det S_J.
\end{equation*}
Put $I=\{i_1,\dots,i_n\}$ and $J=\{j_1,\dots,j_n\}$.
Then there exist partitions $\lambda=(\lambda_1,\dots,\lambda_n)$ and $\mu=(\mu_1,\dots,\mu_n)$ such that $\lambda_1,\mu_1\leq N-n$
and we can write $N-i_r=\lambda_r+n-r$ and $N-j_r=\mu_r+n-r$  for $r=1,\dots,n$.
Then it is easy to see that $\det T_I=\det\left(x_i^{\lambda_j+n-j}\right)=\Delta(x)s_{\lambda}(x)$ and $\det S_J=\det\left(y_i^{\mu_j+n-j}\right)=\Delta(y)s_{\mu}(y)$.
As before we put $k_r=\lambda_r+n-r$ and $\ell_r=\mu_r+n-r$ for $r=1,\dots,n$.
Then, by Lemma~\ref{lemma:coeff},
we obtain
\begin{eqnarray*}
\det A^I_J
=\det\left[\frac{(-q^2;q^2)_{k_i}}{(q^2;q^2)_{k_i}}\frac{(-q^2;q^2)_{\ell_j}}{(q^2;q^2)_{\ell_j}}\frac2{q^{k_i-\ell_j}+q^{\ell_j-k_i}}\right]
=q^{|\lambda-\mu|+|\mu-\lambda|}J_{\lambda\mu}(q^2).
\end{eqnarray*}
On the other hand,
by Lemma~\ref{lemma:two-vars},
we have
\begin{eqnarray*}
\lim_{N\rightarrow\infty}\det{}^tTAS
=\det\left[\frac{(-qx_i;q^2)_{\infty}}{(qx_i;q^2)_{\infty}}\frac{(-qy_j;q^2)_{\infty}}{(qy_j;q^2)_{\infty}}\frac{1}{1-x_iy_j}\right]_{i,j=1,\dots,n}.
\end{eqnarray*}
Thus \thetag{\ref{eqn:kawanaka}} is an immediate consequence of \thetag{\ref{Cauchy2}}.
This proves the theorem.
$\Box$

%% file: appendix.tex
%
%
%
%
%
%

\section{Appendix: A variant of the Sundquist formula}

We give here some variant (both a statement and a proof) of the
Sundquist formula \cite{Su2}. Indeed, we establish the following 
Theorem \ref{Sundquist}. 
Although the initial proof of the theorem 
was made by employing the basic identity in \S2, it did not 
use directly 
the minor summation formula  
and was also, in fact, complicated. Thus we decided to treat this  
in the Appendix. The proof presented here is the one 
followed by the suggestion of the referee. 

\begin{thm}\label{Sundquist}
It holds that 
\begin{align}
&\Pf\left(\frac{y_i-y_j}{a+b(x_i+x_j)+cx_ix_j}\right)_{1\le i,j\le2n}
\hskip-3mm\,\,\times\prod_{1\le i<j\le2n}\{a+b(x_i+x_j)+cx_ix_j\}
\nonumber\\
&=(ac-b^2)^{\binom{n}2}
\sum_{{I\subseteq[2n]}\atop{\sharp I=n}}
(-1)^{|I|-\binom{n+1}2}y_{I}
\Delta_{I}(x)\Delta_{\overline I}(x)J_{I}(x)J_{\overline I}(x)
\nonumber\\
&=(ac-b^2)^{\binom{n}2}
\sum_{{I\subseteq[2n]}\atop{{\sharp I=n}\atop{i_1<j_1}}}
(-1)^{|I|-\binom{n+1}2}(y_{I}+(-1)^ny_{\overline I})\Delta_{I}(x)
\Delta_{\overline I}(x)J_{I}(x)J_{\overline I}(x),
\end{align}
where the sum runs over all $n$-elements subset $I=\{i_1<\cdots<i_n\}$ of
$[2n]=\{1,2,\ldots, 2n\}$ such that $i_1<j_1$ and $|I|=i_1+\cdots+i_n$.
Moreover $\overline I=\{j_1<\cdots<j_n\}$ is the complementary subset of
$I$ in $[2n]$ and
\begin{align*}
&\Delta_{I}(x)=\prod_{{i,j\in I}\atop{i<j}} (x_i-x_j),\\
&J_{I}(x)=J_{I}(x; a,b,c)
=\prod_{{i,j\in I}\atop{i<j}} \{a+b(x_i+x_j)+cx_ix_j\},\\
&y_I=\prod_{i\in I}y_i.
\end{align*}
In particular, if the relation $ac=b^2$ holds then 
\begin{align*}
\Pf\left(\frac{y_i-y_j}{a+b(x_i+x_j)+cx_ix_j}\right)_{1\le i,j\le2n}=0.
\end{align*}
\end{thm}
%
%
\begin{ex}
In the case of $n=2$, if we put $a=c=1$ and $b=0$ then 
the theorem above reads
\begin{align*}
&\Pf\left[\frac{y_{i}-y_{j}}{1+x_{i}x_{j}}\right]_{1\leq i<j\leq4}
\times\prod_{1\leq i<j\leq4}(1+x_{i}x_{j})\\
&=(y_{1}y_{2}+y_{3}y_{4})(x_{1}-x_{2})(x_{3}-x_{4})(1+x_{1}x_{2})(1+x_{3}x_{4})\\
&-(y_{1}y_{3}+y_{2}y_{4})(x_{1}-x_{3})(x_{2}-x_{4})(1+x_{1}x_{3})(1+x_{2}x_{4})\\
&+(y_{1}y_{4}+y_{2}y_{3})(x_{1}-x_{4})(x_{2}-x_{3})(1+x_{1}x_{4})(1+x_{2}x_{3})
\end{align*}
\end{ex}
\vskip .1in
%
%
\noindent
{\it Proof of Theorem \ref{Sundquist}.}
Since 
\begin{equation*}
a+b(x_i+x_j)+cx_ix_j=(\sqrt cx_i + \frac{b}{\sqrt c})(\sqrt cx_j + \frac{b}{\sqrt c}) + a- \frac{b^2}c
\end{equation*}
it is enough to show the theorem for the case  $a=c=1$ and $b=0$.
 
Moreover, since the second equality follows immediately from the 
first one if one notes the fact that $|I| + |\overline I| \equiv n \mod 2$, 
we give a proof of the first equality. First, we notice the following
\begin{lem}\label{yI}
The coefficient of $y_I$ in the Pfaffian 
$\Pf\left(\frac{y_j-y_i}{1+x_ix_j}\right)$  
is equal to the Pfaffian of the following skew symmetric matrix $T_I$:
The $(i,j)$ entry $(T_I)_{ij}$ of $T_I$ is given by
\begin{align*}
(T_I)_{ij}= \begin{cases}
    1/(1+x_ix_j)  & \quad \text{if}\quad i \in I 
\;\text{and}\; j \in \overline I,\\
    -1/(1+x_ix_j) & \quad \text{if}\quad i 
\not\in \overline I \;\text{and}\; j \in I,\\
    0             & \quad \text{otherwise}.
\end{cases}
\end{align*}
\end{lem}
\noindent
{\it Proof.} Recall the definition of a Pfaffian: 
$$
\Pf\left(\frac{y_j-y_i}{1+x_ix_j}\right)
=\sum_{\sigma} \epsilon(\sigma)\,
\frac{y_{\sigma_1}-y_{\sigma_2}}{1+x_{\sigma_1}x_{\sigma_2}}
\cdots \frac{y_{\sigma_{2n-1}}-y_{\sigma_{2n}}}
{1+x_{\sigma_{2n-1}}x_{\sigma_{2n}}},
$$
where 
the summation is over all partitions $\sigma=\{\{\sigma_{1},\sigma_{2}\}_{<},\hdots,\{\sigma_{2n-1},\sigma_{2n}\}_{<}\}$
of $[2n]$ into $2$-elements blocks,
and $\epsilon(\sigma)=\epsilon(\sigma_{1},\hdots ,\sigma_{2n})$ 
denotes the sign of $\sigma\in \frak S_{2n}$. 
From this expression we immediately see that 
the coefficient of $y_I$ 
is given by
$$
\sum_{\sigma} \epsilon(\sigma)
\prod_{\sigma_{2k-1}\in I\atop{\sigma_{2k}\in \overline I}}
\frac1{1+x_{\sigma_{2k-1}}x_{\sigma_{2k}}}
\prod_{\sigma_{2k-1}\in \overline I\atop{\sigma_{2k}\in I}}
\frac{-1}{1+x_{\sigma_{2k-1}}x_{\sigma_{2k}}}.
$$
Note here that when $\sigma$ is subject to 
either the conditions  
$\sigma_{2k-1}\in I, \sigma_{2k}\in I$ or 
$\sigma_{2k-1}\in \overline I, \sigma_{2k}\in \overline I$, 
the corresponding term disappears in the sum. 
Hence the assertion follows easily. $\square$

\smallskip

By this lemma, in order to prove the theorem (for $a=c=1,\, b=0$),
 it suffices to show that  
\begin{align}\label{TI-pfaffian}
\Pf(T_I)=(-1)^{|I|-\binom{n+1}2}
\Delta_{I}(x)\Delta_{\overline I}(x)
\prod_{i\in I, j\in \overline I}\frac1{1+x_ix_j}.
\end{align}
This Pfaffian $\Pf(T_I)$ is computed by using a relation between 
Pfaffians of special type and determinants, and the Cauchy 
determinant formula as follows: 
Recalling $I=\{i_1<i_2<\hdots<i_n\}\subseteq[2n]$ and 
$\overline I=\{j_1<j_2<\hdots<j_n\}\subseteq[2n]$, we first  
notice that
$$
\Pf(T_I)=(-1)^{|I|-\binom{n+1}2}
\Pf\begin{pmatrix} 0 & X_I\\ -X_I&0 \end{pmatrix},
$$
where the $n\times n$ matrix $X_I$ is determined by 
$(X_I)_{k\ell}
=\frac1{1+x_{i_k}x_{j_\ell}}\;(i_k\in I,\, j_\ell\in \overline I)$ 
because the number of the column-row changes for obtaining 
$\Pf\begin{pmatrix} 0 & X_I\\ -X_I&0 \end{pmatrix}$ 
from $\Pf(T_I)$ equals 
$(i_1-1)+(i_2-2)+\cdots+(i_n-n)
=|I|-\frac12n(n+1)=|I|-\binom{n+1}2$. 
Moreover, since
$$
\Pf\begin{pmatrix} 0 & X\\ -X&0 \end{pmatrix}
=(-1)^{\binom{n}2}\det(X)
$$
for any $n\times n$ matrix $X$, the Cauchy determinant formula
claims that 
\begin{align*}
\Pf\begin{pmatrix} 0 & X_I\\ -X_I&0 \end{pmatrix}
&=(-1)^{\binom{n}2}
\Delta_{I}(-x)\Delta_{\overline I}(x)
\prod_{i\in I, j\in \overline I}\frac1{1-(-x_i)x_j}\\
&=\Delta_{I}(x)\Delta_{\overline I}(x)
\prod_{i\in I, j\in \overline I}\frac1{1+x_ix_j},
\end{align*}
whence the equation \eqref{TI-pfaffian} follows. Multiplying 
the factor $\prod_{1\leq i,j\leq 2n}(1+x_ix_j)$ 
to the both sides of \eqref{TI-pfaffian}, we 
obtain the desired identity. This proves the 
theorem. $\square$

\medskip


As a corollary of this theorem we obtain the Sundquist
identity \cite{Su2}.  
The Sundquist identity is a two-variable 
generalization of $\Pf(\frac{x_j-x_i}{1-tx_ix_j})$ and it is 
considered as a Pfaffian version of  Cauchy determinant formula,
whose evaluation is given by \cite{Ste} (see also Lemma 8 in \cite{IW1}): 
\begin{equation*}
\Pf(\frac{x_j-x_i}{1-tx_ix_j})_{1\leq i<j\leq 2n} 
= t^{n(n-1)}\frac{\prod_{1\leq i<j\leq 2n}(x_j-x_i)}{\prod_{1\leq i<j\leq
2n}(1-tx_ix_j)}. 
\end{equation*}

\begin{cor} (Sundquist)
\begin{equation*}
\Pf\left(\frac{y_i-y_j}{1+x_ix_j}\right)_{1\le i,j\le2n}
\times\prod_{1\le i<j\le2n}(1+x_ix_j)
=\sum_{\lambda,\mu}
a_{\lambda+\delta_n,\mu+\delta_n}(x,y),
\end{equation*}
where the sums runs over pairs of partitions
\begin{equation*}
\lambda=(\alpha_1,\cdots,\alpha_p|\alpha_1+1,\cdots,\alpha_p+1),
\mu=(\beta_1,\cdots,\beta_p|\beta_1+1,\cdots,\beta_p+1)
\end{equation*}
in Frobenius notation with $\alpha_1,\beta_1<n-1$.
Also, for $\alpha$ and $\beta$ partitions (compositions, in general)  of
length $n$, we put 
\begin{equation*}
a_{\alpha,\beta}(x,y)
=\sum_{\sigma\in \frak
S_{2n}}\epsilon(\sigma)\sigma(x_1^{\alpha_1}y_1\cdots
x_n^{\alpha_n}y_nx_{n+1}^{\beta_1}\cdots x_{2n}^{\beta_n}),
\end{equation*}
where $\sigma\in \frak S_{2n}$ acts on each of two sets of variables $\{x_1,\cdots,x_n\}$ and $\{y_1,\cdots,y_n\}$ by permuting indices,
and $\delta_n=(n-1,n-2,\cdots,0)$.
$\Box$
\end{cor}
We have already given in \cite{IW5} 
a way of reduction of this corollary 
 from Theorem \ref{Sundquist} by using the expansion 
\begin{equation*}
\prod_{1\le i<j\le n}(1+x_ix_j)
=\sum_{\lambda=(\alpha_1,\cdots,\alpha_p|\alpha_1+1,\cdots,\alpha_p+1)}
s_{\lambda}(x_1,\cdots,x_n),
\end{equation*}
where $s_\lambda=s_\lambda(x_1,\cdots,x_n)$ are the Schur functions, 
so we omit the proof.

\smallskip

Suppose $n$ is even. If we put $y_i=1$ for all $1\leq i\leq
2n$ in Theorem \ref{Sundquist}, it is immediate to see the
\begin{cor} 
\begin{equation*}
\sum_{{I\subseteq[2n]}\atop{{\sharp I=n}\atop{i_1<j_1}}}
(-1)^{|I|}\Delta_{I}(x)\Delta_{I}(x)J_{I}(x)J_{\overline I}(x)=0
\end{equation*}
holds for even $n$.
$\Box$
\end{cor}

%
%
\begin{ex} When $n=2$ we have 
\begin{align*}
&(x_{1}-x_{2})(x_{3}-x_{4})(1+x_{1}x_{2})(1+x_{3}x_{4})\\
-&(x_{1}-x_{3})(x_{2}-x_{4})(1+x_{1}x_{3})(1+x_{2}x_{4})\\
+&(x_{1}-x_{4})(x_{2}-x_{3})(1+x_{1}x_{4})(1+x_{2}x_{3})=0.
\end{align*}
\end{ex}
\vskip .1in

\begin{rem}
It is naturally thought the formula as the identity
of two variables relevant to a $A_n$-type root system. 
It would be interesting to establish the $B_n, C_n, D_n$-analogues
of Theorem \ref{Sundquist} like in \cite{IOW} for the generalization 
of the Littlewood
formulas to the classical groups. 
\end{rem}

\bigskip
\par\noindent
{\bf Acknowledgements.}
The authors would like to thank the referee  
for his/her many valuable comments and suggestions.

%% file: lmsf.bbl
\begin{thebibliography}{99}


\bibitem{B1}
D.M.~Bressoud,
{\em Proofs and Confirmations},
Cambridge U.P.

\bibitem{B2}
D.M.~Bressoud,
{\em Identities for Schur functions and plane partitions},
Ramanujan J. {\bf 4} (2000) 69-80.

\bibitem{DW}
A.~Dress and W.~Wenzel,
{\em A simple proof of an identity concerning pfaffians of skew symmetric matrices},
Adv. Math.
{\bf 112} (1995), 120--134.


\bibitem{GNT}
A.~Gyoja, K.~Nishiyama and K.~Taniguchi,
{\em Kawanaka invariants for representations of Weyl groups},
J. Alg., 225 (2000), 842 - 871.


\bibitem{GV}
I.~Gessel and G.~Viennot,
{\em Determinants, Paths, and Plane Partitions},
preprint (1989).



\bibitem{Hi}
R.~Hirota,
{\em Mathematical aspect of the soliton theory from a direct methods point of view},
Iwanami Shoten
(1992), in Japanese.



\bibitem{IOW}
M.~Ishikawa, S.~Okada and M.~Wakayama,
{\em Applications of minor summation formulas I,
Littlewood's formulas},
J. Alg.
{\bf  183} (1996), 193--216.



\bibitem{IW1}
M.~Ishikawa and M.~Wakayama,
{\em Minor summation formula of Pfaffians},
Linear and Multilinear Alg.
{\bf  39} (1995), 285-305.




\bibitem{IW2}
M.~Ishikawa and M.~Wakayama,
{\em Minor summation formula of Pfaffians and Schur functions identities},
Proc.Japan Acad., Ser.A
{\bf  71} (1995), 54--57.



\bibitem{IW3}
M.~Ishikawa and M.~Wakayama,
{\em New Schur function series},
J. Alg, 
{\bf  208} (1998), 480--525.



\bibitem{IW4}
M.~Ishikawa and M.~Wakayama,
{\em Applications of minor summation formulas II,
Pfaffians of Schur polynomials},
J. Combin. Th. Ser.~A
{\bf  88} (1999), 136--157.



\bibitem{IW5}
M.~Ishikawa and M.~Wakayama,
{\em Minor summation formula of Pfaffians, Survey and a new identity},
Adv. Stud. Pure Math.
{\bf  28} (2000), 133--142.


\bibitem{JZ1}
F.~Jouhet and J.~Zeng,
{\em Some New Identities for Schur Functions},
Adv. in Applied Math.
{\bf  27} (2001), 493--509.


\bibitem{JZ2}
F.~Jouhet and J.~Zeng,
{\em New Identities of Hall-Littlewood Polynomials and Applications},
to appear in Ramanujan J.


\bibitem{Ka1}
N.~Kawanaka,
{\em A $q$-series identity involving Schur functions and related topics},
Osaka J. Math.
{\bf  36} (1999), 157--176.

\bibitem{Ka2}
N.~Kawanaka,
{\em A $q$-Cauchy identity involving Schur functions and imprimitive complex reflection groups},
Osaka J. Math.
{\bf  38} (2001), 775--810.



\bibitem{KW}
K.~Kinoshita and M.~Wakayama,
{\em Explicit Capelli identities for skew symmetric matrices},
Proc. Edinburgh Math. Soc. 
{\bf 45} (2002), 449--465.


\bibitem{Kn}
D.~Knuth,
{\em Overlapping pfaffians},
Electronic J. of Combi.
{\bf  3} (1996), 151--163.




\bibitem{Kr1}
C.~Krattenthaler,
{\em Determinant identities and a generalization of the number of totally symmetric self-complementary plane partitions},
Electron. J. Combin.
{\bf 4} (1) (1997), $\sharp$R27, 62.

\bibitem{Kr2}
C.~Krattenthaler,
{\em Advanced Determinant Calculus},
Seminaire Lotharingien Combin.
{\bf 42} ("The Andrews Festschrift") (1999), Article B42q, 67.





\bibitem{La}
A.~Lascoux,
{\em Littlewood's formulas for characters of orthogonal and symplectic groups},
Algebraics Combinatorics and Quantum Groups (2003).



\bibitem{LT1}
J.-G. Luque and J.-Y. Thibon,
{\em Hankel hyperdeterminants and Selberg integrals}
J. Phys. A: Math. Gen. 36 (16 May 2003) 5267-5292.

\bibitem{LT2}
J.-G. Luque and J.-Y. Thibon,
{\em Pfaffian and hafnian identities in shuffle algebras},
Adv. Appl. Math. {\bf 29} (2002), 620--646.


\bibitem{Ma}
I.~G.~Macdonald,
{\em Symmetric Functions and Hall Polynomials,
 2nd Edition},
Oxford University Press,
(1995).


\bibitem{MRR1}
W.H.~Mills, D.P.~Robbins and H.~Rumsey, Jr.,
{\em Self-Complementary Totally Symmetric Plane Partitions},
J. Combin. Th. Ser.~A
{\bf  42} (1986), 277--292.


\bibitem{O1}
S.~Okada,
{\em On the generating functions for certain classes of plane partitions},
J. Combin. Th. Ser.~A
{\bf  51} (1989), 1--23.



\bibitem{O2}
S.~Okada,
{\em Applications of minor-summation formulas to rectangular-shaped representations  of classical groups},
J. Alg.
{\bf  205} (1998), 337--367.


\bibitem{RR}
D.P.~Robbins and H.~Rumsey Jr.,
{\em Determinants and Alternating Sign Matrices},
Adv. Math.
{\bf 62} (1986), 169--184.


\bibitem{S}
I.~Schur,
Aufgabe 569, 
Arch. Math. Phys. (3) {\bf 27} (1918), 163.



\bibitem{Sta}
R.~P.~Stanley,
{\em Enumerative Combinatorics, Volume I, II},
Cambridge University Press,
(1997), (1999).


\bibitem{Ste}
J.~Stembridge,
{\em  Nonintersecting paths, pfaffians and plane partitions},
Adv. Math.
{\bf  83} (1990), 96--131.






\bibitem{Su2}
T.~Sundquist
{\em  Two variable Pfaffian identities and symmetric functions}
J. Alg. Combin.
{\bf  5} (1996), 135--148.



\bibitem{TW1}
C.A.~Tracy and H.~Widom,
{\em A Limit Theorem for Shifted Schur Measures},
Duke Math. J. {\bf 123} (2004) 178--208.











\end{thebibliography}
